\documentclass[12pt]{amsart}
\usepackage{amssymb}

\def\thesisform{0}

\def\withfig{0}

\ifnum\withfig=1
\usepackage{psfig}
\fi

\newtheorem{thm}{Theorem}[section]
\newtheorem{lem}{Lemma}[section]
\newtheorem{cor}{Corollary}[section]
\newtheorem{conj}{Conjecture}[section]
\newtheorem{prop}{Proposition}[section]

\theoremstyle{remark}

\newtheorem{rem}{Remark}[section]
\newtheorem{claim}{Claim}[section]

\numberwithin{equation}{section}

\allowdisplaybreaks[4]

\ifnum\thesisform=1
\fi

\begin{document}

\ifnum\thesisform=1
\begin{titlepage}
\begin{center}

\vspace{3in}
 
{\Large {\bf Rational Curves on K3 Surfaces}}

\vspace{.3in}

{\large A thesis presented

\vspace{.2in}

by

\vspace{.2in}

Xi Chen

\vspace{.2in}

to

\vspace{.2in}

The Department of Mathematics

\vspace{.2in}

in partial fulfillment of the requirements

%\vspace{.1in}

for the degree of

%\vspace{.1in}

Doctor of Philosophy

%\vspace{.1in}

in the subject of

%\vspace{.1in}

Mathematics

\vspace{.3in}

Harvard University

%\vspace{.1in}

Cambridge, Massachusetts

%\vspace{.1in}

March 1997}
\end{center}
\newpage
\end{titlepage}

\pagenumbering{roman}
\thispagestyle{empty}

\begin{center}
\quad
\vspace{3in}

{\copyright 1997 by Xi Chen

\vspace{.1in}

All rights reserved.}
\end{center}

\newpage

\baselineskip=24pt
\fi

\pagenumbering{arabic}

% macros

\newcommand{\claimref}[1]{Claim \ref{#1}}
\newcommand{\thmref}[1]{Theorem \ref{#1}}
\newcommand{\propref}[1]{Proposition \ref{#1}}
\newcommand{\lemref}[1]{Lemma \ref{#1}}
\newcommand{\coref}[1]{Corollary \ref{#1}}
\newcommand{\remref}[1]{Remark \ref{#1}}
\newcommand{\conjref}[1]{Conjecture \ref{#1}}
\newcommand{\secref}[1]{Sec. \ref{#1}}
\newcommand{\ssecref}[1]{\ref{#1}}
\newcommand{\sssecref}[1]{\ref{#1}}

\def \d#1{\displaystyle{#1}}
\def \mult{\mathop{\mathrm{mult}}\nolimits}
\def \rank{\mathop{\mathrm{rank}}\nolimits}
\def \codim{\mathop{\mathrm{codim}}\nolimits}
\def \Ord{\mathop{\mathrm{Ord}}\nolimits}
\def \Var{\mathop{\mathrm{Var}}\nolimits}
\def \Ext{\mathop{\mathrm{Ext}}\nolimits}
\def \EQ{\Leftrightarrow}
\def \Pic{\mathop{\mathrm{Pic}}\nolimits}
\def \Spec{\mathop{\mathrm{Spec}}\nolimits}
\def \mapright#1{\smash{\mathop{\longrightarrow}\limits^{#1}}}
\def \mapleft#1{\smash{\mathop{\longleftarrow}\limits^{#1}}}
\def \mapdown#1{\Big\downarrow\rlap{$\vcenter{\hbox{$\scriptstyle#1$}}$}}
\def \smapdown#1{\downarrow\rlap{$\vcenter{\hbox{$\scriptstyle#1$}}$}}
\def \A{{\mathbb A}}
\def \I{{\mathcal I}}
\def \J{{\mathcal J}}
\def \CO{{\mathcal O}}
\def \C{{\mathcal C}}
\def \BC{{\mathbb C}}
\def \m{{\mathcal M}}
\def \H{{\mathcal H}}
\def \S{{\mathcal S}}
\def \Z{{\mathcal Z}}
\def \BZ{{\mathbb Z}}
\def \Y{{\mathcal Y}}
\def \T{{\mathcal T}}
\def \P{{\mathbb P}}
\def \G{{\mathbb G}}
\def \F{{\mathbb F}}
\def \BR{{\mathbb R}}
\def \Hom{{\mathrm{Hom}}}
\def \nilrad{{\mathrm{nilrad}}}
\def \Supp{{\mathrm{Supp}}}
\def \Ann{{\mathrm{Ann}}}
\def \closure#1{\overline{#1}}
\def \EQ{\Leftrightarrow}
\def \imply{\Rightarrow}
\def \isom{\cong}
\def \embed{\hookrightarrow}
\def \tensor{\otimes}
\def \wt#1{{\widetilde{#1}}}
\def \Jac{{\mathrm{Jac}}}

\title{Rational Curves on K3 surfaces}
\author{Xi Chen}
\address{UCLA Department of Mathematics\\
6363 Math Sciences\\
Box 951555\\
Los Angeles, CA 90095-1555}
\email{xchen@math.ucla.edu}
\date{\today}
\maketitle

\section{Introduction}

The classification theory of algebraic surfaces shows there are at most
countably many rational curves on a K3 surface. The first question we may
ask is whether there are any rational curves at all.
The existence of rational curves on a general K3 surface was established in
\cite{M-M}. A generalization was made by S. Nakatani as follows.

\begin{thm}[Nakatani]
Let $F_m$ be the moduli of pairs $(S,L)$ of a K3 surface $S$
and a non-divisible ample divisor $L$ on it such that $L^2=2m-2$
modulo obvious isomorphisms. If $m$ is odd then, for a sufficiently general
$(S,L) \in F_m$, $|(k^2+1)/2L|$ has an irreducible
rational curve for every odd number $k$.
\end{thm}

In \secref{S:EXIST} we will extend the existence of irreducible
rational curves to every complete linear series on a general K3 surface,
i.e.,

\begin{thm}\label{T:EXIST}
For any integers $n\ge 3$ and $d > 0$,
the linear system $|\CO_S(d)|$ on a general K3 surface $S$ in ${\mathbb
P}^n$ contains an irreducible nodal rational curve.
\end{thm}

It must be mentioned that this is a folklore theorem known to several
people. But no complete proof has appeared in literature yet.

The next natural question following the existence problem is how many
irreducible rational curves there are in $|\CO(d)|$ on a general K3 surface
in ${\mathbb P}^n$. The number for $d=1$ has been
successfully calculated in \cite{Y-Z}. They give the following
remarkable formula
\begin{equation}\label{E:YAU-ZASLOW}
\sum_{g=1}^\infty n(g) q^g = \frac{q}{\Delta(q)}
\end{equation}
where $\Delta(q) = q\prod_{n=1}^\infty (1-q^n)^{24}$ is the well-known
modular form of weight $12$ and $n(g)$ is the nominated number of rational
curves in $|\CO(1)|$ on a general K3 surface in $\P^g$ for $g\ge 3$. 
More precisely, $n(g)$ is the
sum of the Euler characteristics of the compactified Jacobians
of all rational curves in
$|\CO(1)|$ (for a detailed exposition, see \cite{B}).
Since the compactified Jacobian of a rational curve with singularities other
than nodes is not very well understood, we only know this sum equals the
number of rational curves in $|\CO(1)|$ on a K3 surface in the case that
all these rational curves are nodal. Hence
the only gap left in this enumeration problem
is the hypothesis that all rational curves in $|\CO(1)|$
on a general K3 surface are nodal, namely, the following conjecture,

\begin{conj}\label{T:NOD}
For $n\ge 3$,
all rational curves in the linear system $|\CO_S(1)|$ on a general K3 surface
$S$ in ${\mathbb P}^n$ are nodal.
\end{conj}

A proof of \conjref{T:NOD} has not been completely worked out at the time
this paper is written. We will show the readers our approach and progress
made towards this problem.

Basically, we study rational curves on K3 surfaces here by specializing K3
surfaces. We will
degenerate general K3 surfaces to some ``special'' ones and study
limits of rational curves on these special K3 surfaces. By examining these
``limiting rational curves'', hopefully we can say something about rational
curves on a general K3 surface.

To be specific, we will degenerate a general K3 surface to a trigonal K3
surface (see \secref{S:TRIG1} for definition) in order to show
\conjref{T:NOD}. A main theorem (\thmref{T:NODB}) will be proved, which
enable us to convert \conjref{T:NOD} into some similar statements
(\conjref{T:NODA} and \ref{T:NODC}) concerning rational curves on a
trigonal K3 surface. By a corollary of \thmref{T:NODB}, we see that
\conjref{T:NOD} is true for $n\le 9$ and $n=11$ and hence
justify the Yau-Zaslow's counting formula \eqref{E:YAU-ZASLOW} for
$g\le 9$ and $g = 11$. This statement is also proved independently by Kang Zuo.

Z. Ran is working on a related problem concerning 
curves of any genus on a quartic
surface. He has obtained

\begin{thm}[Ran]
Rational curves of any degree on a quartic surface have transitive
monodromy. Namely, if we let $W_d$ be the correspondence $(C, S)$ that $C\in
|\CO_S(d)|$ is a rational curve on a quartic surface $S\in
|\CO_{\P^3}(4)|$, then $W_d$ is irreducible for any $d>0$.
\end{thm}

With this result in mind, it is reasonable to conjecture that

\begin{conj}\label{conj:1}
Rational curves in $|\CO_S(d)|$ on a primitive K3 surface $S\subset \P^n$ have
transitive monodromy for any $d>0$.
\end{conj}

Notice that since we already have \thmref{T:EXIST},
\conjref{conj:1} implies that every rational curves on a
general K3 surface is nodal, which is much stronger than \conjref{T:NOD}.

\medskip\noindent{\bf Conventions.}

\begin{enumerate}
\item Throughout the paper, we will work exclusively over 
${\mathbb C}$.
\item We will only concern ourselves with the primitive K3
surfaces here. Hence from time to time we will
simply call a general primitive K3 surface in $\P^n$
a general K3 surface in $\P^n$.
Hopefully no confusion would arise from this abuse of terminology.
\item Since we are working over ${\mathbb C}$, we will use analytic geometry
whenever possible. Hence we will use analytic neighborhoods of points
instead of Zariski open neighorhoods in most cases,
while you may always replace them by formal or etale neighborhoods.
\item A double curve singularity $x^2 - y^{n+1} = 0$ is also called an
$A_n$ singularity under the A-D-E classification of simple singularities. 
Here we allow $n=0$ in $A_n$ which simply refers to a smooth point.
\end{enumerate}

\medskip\noindent{\bf Acknowlegments.}

I would like to thank Joe Harris, my thesis advisor, for suggesting the
problem and help me throughout this paper. I also benefited greatly from the
discussions with Eric Zaslow, Ciro Cilibeto and Tony Pantev. 

I also thank
Angelo Vistoli for reading through the first draft of this paper and
providing several valuable suggestions.

\section{Preliminaries}

This section is a miscellaneous collection of
theorems and results which,
though needed in our proof of \thmref{T:EXIST} and later our attempt to the
proof of \conjref{T:NOD}, do not fit
very well there. We place
them here so that they will not disrupt our main course of discussion
later. Readers are suggested to skip this section and only come back when
the results stated in this section are referred.

\subsection{Nodal reduction of a family of curves}\label{s:nodred}

Let $\Upsilon\to T$ be an irreducible family of curves over disk $T$. We
want to introduce a common construction on $\Upsilon$ which will be used
througout the paper. 

Let $\wt{\Upsilon}$ be the normalization of the surface $\Upsilon$.
After an apporiate base
changes $\pi: \Delta\to T$ and necessary blow-ups, 
we arrive
at a family $\Upsilon^v\to \Delta$ with the diagram
\begin{equation*}
\begin{array}{ccccc}
\Upsilon^v & \mapright{v} & \wt{\Upsilon}  & \mapright{} & \Upsilon \\
\downarrow &         & \downarrow & \swarrow \\
\Delta & \mapright{\pi} & T
\end{array}
\end{equation*}
and the properties
\begin{enumerate}
\item $\Upsilon^v$ is smooth;
\item the central fiber $\Upsilon_0^v$
of $\Upsilon^v\to\Delta$ only has nodes as its singularities;
\item $\Upsilon^v\to \Delta$ is minimal with respect to these properties.
\end{enumerate}
Obviously, the minimality condition means that there are no contractible $-1$
rational curves on $\Upsilon_0^v$, i.e., there are no rational components of
$\Upsilon_0^v$ mapped constantly to $\wt{\Upsilon}$ by $v$ and meeting the
rest of $\Upsilon_0^v$ only at one point.

Depending on our needs, we may further blow down the contractible
$-2$ rational curves
of $v: \Upsilon^v\to \wt{\Upsilon}$, namely, the rational components of
$\Upsilon_0^v$ mapped constantly to $\wt{\Upsilon}$ by $v$ and meeting the
rest of $\Upsilon_0^v$ at exactly two points. In this case, we have to drop the
smoothness of the total family and replace it by
\begin{enumerate}
\item[$1'$.] $\Upsilon^v$ has only isolated double points as its singularities.
\end{enumerate}

We may also ``mark'' $\wt{\Upsilon} \to T$ with $n$ different sections
$s_i: T\to \wt{\Upsilon}$ for $i=1,2,...,n$
(typically they come from the singular locus of
$\Upsilon \to T$). We can ask $v: \Upsilon^v\to \wt{\Upsilon}$
to separate the $n$ sections on the central fiber, namely,
\begin{enumerate}
\item[$2\frac{1}{2}$.] there are $n$ sections $s_i^v : \Delta\to
\Upsilon^v$ for $i=1,2,...,n$
such that $v\circ s_i^v = s_i\circ \pi$, $\Upsilon_0^v$ is smooth at point
$s_i^v(\Delta) \cap \Upsilon_0^v$
and $s_i^v(\Delta)\cap \Upsilon_0^v \ne s_j^v(\Delta)\cap
\Upsilon_0^v$ for $i\ne j$.
\end{enumerate}

We will call $\Upsilon^v \to \Delta$ (or simply $\Upsilon^v$) the nodal
reduction of the family $\Upsilon\to T$. We will not distinguish the two
constructions with the smoothness condition 1 or the weaker $1'$ since
in our application either the
difference is unessential or it can be easily
told by the context which one is in use.

\subsection{Deformation of a nonreduced singularity $x^m y^n = 0$}

The following lemma is a merely easy observation 
but it will be very useful in our
furture discussion.

\begin{lem}\label{L:XmYn}
Let $S$ be a family of curves over disk $T$ with irreducible general
fibers, $p$ be a point on the central
fiber $S_0$ of $S$ and $V$ be an analytic (formal or etale) neighborhood of
$p$ on $S$. Suppose that $V$ can be embedded into $\A^2 \times T$ and the
central fiber $V_0$ is correspondingly given by $x^m y^n = 0$ (let
$\A^2 \times T$ be parameterized by $(x, y, t)$). Let $S^v$ be a family of
curves over $T$ and $\pi: S^v\to S$ be a generically 1-1 morphism preserving
base $T$. Suppose that $\pi^{-1}(V)$ consists of $r$ disjoint
irreducible components
$V^1, V^2, ..., V^r$ and the morphism $\pi: V_0^i \to V_0 = \{x^m y^n = 0\}$
on the central fiber
factors through $\{ x^{m_i} y^{n_i} = 0 \}$ where $m_i$ and $n_i$ is minimal
with respect to this property. Then $V$ has $r$ irreducible components and is
given by $\prod_{i=1}^r (x^{m_i} y^{n_i} + O(t)) = 0$.
\end{lem}

\begin{proof}
There is not much to prove. Let $V$ be given by $f(x, y, t) = 0$
in $\A^2\times T$. Since $\pi$
is generically 1-1, $V$ consists of $r$ different irreducible components,
$\pi(V^i)$, for $i=1,2,..., r$. And since $\pi(V^i)$ is obviously given by
$x^{m_i} y^{n_i} + O(t) = 0$, $f(x, y, t) = 
\prod_{i=1}^r (x^{m_i} y^{n_i} + O(t))$.
\end{proof}

\subsection{Deformation of curves on the surface $xy = 0$}

From time to time, we will deal with a family of curves lying on a family
of surfaces whose central fiber consists of two smooth surfaces meeting
transversely along a curve. To set it up, let $X$ be a one-parameter family
of surfaces over disk $T$ with central fiber $X_0 = Q_1\cup Q_2$ where
$Q_1$ and $Q_2$ are two smooth surfaces meeting
transversely along a curve $E=Q_1\cap Q_2$. Let $p$ be a point on
$E$. Since all the theorems stated below are local statements about
families of curves on $X$ at $p$, we
may embed $X$ into ${\mathbb A}^3\times T$ with coordinates $(x, y, z,
t)$ such that $p$ is the origin
$x=y=z=t=0$, $E$ is the line $x=y=t=0$, and $Q_1$ and $Q_2$ be given by
$x=t=0$ and $y=t=0$, respectively. Under these settings, we have

\begin{lem}\label{L:XY=T}
With the setup above, further assume 
that $X$ is locally given by $xy=t^\alpha$ for some integer $\alpha > 0$
at point $p$. Let $S$ be a family of
curves over $T$ and $\pi: S\to X$ be a proper
morphism preserving the base $T$.
Suppose
that there is a curve $C_1$ lying on the central fiber $S_0$ of $S$ mapped by
$\pi$ nonconstantly to a curve on the surface $x=t=0$ passing through $p$.
Then there is correspondingly a curve
$C_2$ lying on the same connected component of $S_0$ as $C_1$
which is mapped by $\pi$ nonconstantly to a curve on the surface
$y=t=0$ passing through $p$.
\end{lem}

\begin{thm}\label{T2}
With the setup above, further assume that $X$ is smooth at $p$. 
Let $L$ be a line bundle on $X$ and $\sigma\subset |L|$ be a linear series
of $L$.

Let $C = C_1\cup C_2$ be a curve on $X_0$ cut out by
an element of $\sigma$, where
$C_i\subset Q_i$ meets $E$ at $p$ with
multiplicity $m\ge 2$ and is smooth at $p$, for $i=1,2$.
Suppose the linear series $\sigma$ generates $m-2$ jets at $p$ on $E$, by
that we mean, the natural map
$$\sigma \subset H^0(X, L)\mapright{} H^0(E, \CO_E/{\mathcal M}_p^{m-1}\tensor
L)\mapright{} 0$$
is surjective, where ${\mathcal M}_p$ is the maximal ideal of $p$ on $E$. 

Let $U$ be an analytic neighborhood of $p$ where $X$ is smooth and $C$ is
smooth outside of $p$, $Y\subset \sigma\times (T-\{0\})$ be defined by
$$Y=\{(s, t): s\in \sigma, t\ne 0, \ {\mathrm{the\ curve}}\
\{s=0\}\cap X_t \ {\mathrm{has}}\ m-1 \ {\mathrm{nodes\ in}}\ U\},$$
and $\closure{Y}$ be the closure of $Y$ in $\sigma\times T$. If $s_0\in
\sigma$ cuts out $C$ on $X_0$,
then
\begin{enumerate}
\item $(s_0, 0)\in \closure{Y}$ and hence $Y$ is nonempty; and $Y$ has
codimension $m-1$ in $\sigma\times T$;
\item the central
fiber of $\closure{Y}\to T$ is nonreduced with multiplicity $m$ in a
neighborhood of $(s_0, 0)$.
\end{enumerate}
\end{thm}

\begin{thm}\label{T3}
With the setup above, 
further assume that $p$ is a rational double point of $X$.
Let $L$ be a line bundle on $X$ and $\sigma\subset |L|$ be a linear series
of $L$.

Let $C = C_1\cup C_2$ be a curve on $X_0$ cut out by
an element of $\sigma$, where
$C_i\subset Q_i$ has 
an ordinary singularity of multiplicity $m>0$ at $p$ and every branch of
$C_i$ at $p$ intersects $E$ at $p$ transversely, for $i=1,2$.
Suppose the linear series $\sigma$ generates $(m-1)$-jets at $p$ on 
$X_0$, by that we mean, the natural map
$$\sigma\subset H^0(X, L)\mapright{} H^0(X_0, \CO_{X_0}/{{\mathcal M}_p}^m\tensor
L)\mapright{} 0$$
is surjective, where ${\mathcal M}_p$ is the maximal ideal of $p$ on $X_0$.

Let $U$ be an analytic neighborhood of $p$ where $X$ and $C$ are
smooth outside of $p$, $Y\subset \sigma\times T-\{0\}$ be defined by
$$Y=\{(s, t): s\in \sigma, t\ne 0, \ {\mathrm{the\ curve}}\
\{s=0\}\cap X_t \ {\mathrm{has}}\ m^2 \ {\mathrm{nodes\ in}}\ U\},$$
and $\closure{Y}$ be the closure of $Y$ in $\sigma\times T$. If
$s_0\in\sigma$ cuts out $C$ on $X_0$,
\begin{enumerate}
\item $(s_0, 0)\in \closure{Y}$ and hence $Y$ is nonempty; and $Y$ has
codimension $m^2$ in $\sigma\times T$; 
\item the central
fiber of $\closure{Y}\to T$ is irreducible and smooth in a neighborhood of
$(s_0, 0)$.
\end{enumerate}
\end{thm}

\begin{proof}[Proof of \lemref{L:XY=T}]
We may assume that $S$ is smooth and $\pi(C_1)$ does not
contain the curve $x=y=t=0$ (otherwise we simply take $C_2 = C_1$).
Let $Q_1$ and $Q_2$ be the surfaces $x=t=0$ and $y=t=0$, respectively.

If $\alpha = 1$, since $X$ is smooth, we have
$$\pi_* (C_1\cdot \pi^* (Q_2)) = \pi_* (C_1)\cdot Q_2 \ne 0.$$ Hence $C_1$
has nonempty intersection with $\pi^{-1}(Q_2)$. The lemma follows.

If $\alpha > 1$, we can resolve the singularities of $X$ by subsequent
blowups and do the induction on $\alpha$.
We can resolve the singularities $X$ as in \cite[Appendix C, p. 39]{G-H} but it
can be done more directly in our case as follows.
Let $\wt{X}$ be the blowup of $X$
along $Q_2$. Then the central fiber of $\wt{X}$ is $(D_1\cup
E\cup D_2)\times \A_z^1$ where $E\isom \P^1$, $D_1\times \A_z^1$ dominates
$Q_1$ and $D_2\times \A_z^1$ dominates $Q_2$. Let $u = t/y$ and $v = y/t$
be the affine coordinates of $E$. Then $D_1$ meets $E$ at point $p_1 = (u =
z = 0)$ where
$\wt{X}$ is given by $uy = t$ and $D_2$ meets $E$ at point $p_2 = (v = z =
0)$ where $\wt{X}$ is given by $xv = t^{\alpha-1}$. Let $\wt{S}$ and
$\wt{C_1}$ be the proper transform of $S$ and $C_1$ and $\wt{\pi}:
\wt{S}\to \wt{X}$. Obviously $\wt{\pi} (\wt{C_1})$ lies on $D_1\times \A_z^1$
and passes through $p_1$. Hence there is a curve $F\subset \wt{S}_0$ lying
on the same connected component as $\wt{C_1}$ and with nonconstant
image $\wt{\pi}(F)\subset E\times \A_z^1$ passing through $p_1$. If
$\wt{\pi}(F)$ contains the curve $u = y = t = 0$ which dominates $x = y=
t= 0$, we are done. If not,
$\wt{\pi}(F)$ must consist of the curve $y = z = t = 0$ and hence passes
through $p_2$. Then by induction hyperthesis, we can find a curve
$C_2'\subset\wt{S}_0$ lying on the same connected component as $F$ and with
nonconstant image $\wt{\pi}(C_2')\subset D_2\times \A_z^1$ passing through
$p_2$. Projecting $C_2'$ to $S_0$, we get $C_2$ as required.
\end{proof}

In both cases described by the \thmref{T2} and \ref{T3},
we will loosely say point $p$ on $C$ can be deformed to
$m-1$ or $m^2$ nodes on a general fiber.
Only the first parts of both theorem are needed for the proof of the
existence theorem, while the second parts will be
useful in our further discussion.

The proof of \thmref{T2} is essentially an application of a theorem of
L. Caporaso and J. Harris \cite{CH2} on deformations of tacnodes.

An $m$-th order tacnode is just another name for the singularity of type
$A_{2m-1}$, which
is analytically equivalent to the origin in the plane
curve given by the equation
$$y(y+x^m) = 0.$$
Without the order specified, a tacnode refers to a second order tacnode.
The versal deformation space of an $m$-th tacnode is then the family $\pi:
\S\to \Delta$, where $\Delta\isom {\mathbb A}^{2m-1}$ with coordinates $(a_{m-2},
...,a_1, a_0, b_{m-1},...,b_0)$, $\S$ is the subscheme of $\Delta\times
{\mathbb A}^2$ given by the equation
$$y^2 +(x^m+a_{m-2}x^{m-2}+...+a_1x+a_0)y + b_{m-1}x^{m-1}+...+b_1x +b_0 =
0$$
and $\pi: \S\to \Delta$ is the projection $\S\subset \Delta\times {\mathbb
A}^2\to \Delta$. Let $\Delta_m$ be the closure of the locus of 
points $(a,b)\in \Delta$ over which the fiber $\S_{a,b}$ has $m$ nodes and
$\Delta_{m-1}$ be the closure of the locus of points $(a,b)$ over which the
fiber $\S_{a,b}$ has $m-1$ nodes.

The theorem of L. Caporaso and J. Harris \cite[Lemma 4.1]{CH2} says

\begin{thm}[Caporaso-Harris]\label{TCH}
Let $m\ge 2$, and let $W\subset\Delta$ be any smooth, $m$-dimensional
subvariety containing $(m-1)$-plane $\Delta_m$, and suppose only that its
tangent plane at the origin is not contained in the hyperplane
$H\subset\Delta$ given by $b_0=0$. Then the intersection
$$W \cap\Delta_{m-1} = \Delta_m \cup \Gamma$$
where $\Gamma$ is a smooth curve having contact of order exactly $m$ with
$\Delta_m$ at the origin.
\end{thm}

\begin{proof}{Proof of \thmref{T2}}
We may assume that $X$ is locally given by $xy = t$ at $p$.
Let $B$ be a neighborhood of $s_0\in \sigma \isom \P^n$, where
$n=\dim\sigma$.
Let $Z\subset B\times X$ be the family of curves over $B\times T$ 
whose fiber over $s\in B$ and $t\in T$ is the curve
$\{s=0\}\cap X_t$.

Let $\S\to\Delta$ be the versal deformation space of an $m$-th tacnode.
Since $C$ has an $m$-th tacnode at $p$, we
have a map $\phi: B\times T\to\Delta$ which induces a local isomorphism
$\psi: Z\mapright{\sim}\S\times_{\Delta} B$ at $p$. 

We claim that the image $\phi(B\times T)$ of $\phi$ is a smooth
$m$-dimensional
subvariety containing $\Delta_m$ and its
tangent plane at the origin is not contained in the hyperplane
given by $b_0=0$.
This can be verified by explicitly writing down the local defining
equations of $Z$ at $p$.
Choose a trivialization of the line bundle $L$ at $p$
and a basis $\{s_0, s_1, ..., s_{m-1}, s_m,..., s_n\}$ of $\sigma$
such that after some scaling of $x,y,z,t$ we have
\begin{eqnarray*}
s_0(x,y,z,t)&=&x+y+z^m+O(t, yz, z^{m+1}),\\
s_i(x,y,z,t)&=&z^{i-1}+O(t, y, z^{m-1}), \ {\mathrm{for}}\ 1\le i\le m-1,\\
s_i(x,y,z,t)&=&O(t, y, z^{m-1}), \ {\mathrm{for}}\ m\le i\le n,
\end{eqnarray*}
where $O(f_1, f_2,...,f_j)$ denotes an element in
the ideal generated by
$\{f_1, f_2,...,f_j\}$.
Let $B$ be parameterized by $(t_1,t_2,...,t_n)$ such that point $(t_1,
t_2,..., t_n)$ represents a section $s=s_0+\sum_{i=1}^n t_i s_i\in
\sigma$. Under these coordinates, $Z$ is locally defined by
$$y\left(y+z^m + \sum_{i=1}^{m-1} t_iz^{i-1} +
O\left(t,yz,z^{m+1},t_iy\big|_{i=1}^n, t_iz^{m-1}\big|_{i=1}^n\right)\right) =
t$$
at $(s_0, p)$.
It is easy to see that
the Kodaira-Spencer map of the family $Z$ at $p$,
i.e., the homomorphism $d\phi$
on the tagent spaces induced by $\phi$,
carries the tagent space of $B\times
T$ at $(s_0,0)$ to the subspace spanned by $\{\partial/\partial a_{m-2},
..., \partial/\partial a_1, \partial/\partial a_0, \partial/\partial b_0\}$
of the tangent space of $\Delta$ at the origin. Hence 
$\phi(B\times T)$ is smooth at the origin with tangent plane $b_{m-1} =
b_{m-2} = ... = b_1 = 0$, which is not contained in the
hyperplane $b_0=0$. And since $(B\times T)\cap \{t=0\}$ is reducible,
we must have $\phi(B\times T)\supset\Delta_m$. 
Then it follows \thmref{TCH} 
that $\phi(B\times T) \cap
\Delta_{m-1} = \Delta_m \cup \Gamma$, where $\Gamma$ is a smooth curve
having contact of order exactly $m$ with $\Delta_m$ at the origin.
Obviously, $Y = \phi^{-1}(\Delta_{m-1}\backslash\Delta_m)$ in $B\times T$. 
Therefore $\closure{Y} = \phi^{-1}(\Gamma)$ and the central fiber of
$\closure{Y} \to T$ is 
$$\closure{Y}\cap \{t=0\} =
\phi^{-1}(\Gamma\cap\Delta_m)$$
and hence must be nonreduced with multiplicity $m$ in $B\times T$.
\end{proof}

The proof of \thmref{T3} involves the study of a non-planary
singularity which, when embedded in $\A^3$, can be put in the form
$$
\left\{
\begin{array}{ccc}
xy &=& 0\\
f(x,y,z) &=& 0,
\end{array}
\right.
$$
where $f(x,y,z)$ is a homogeneous polynomial of order $m$.

\begin{proof}{Proof of \thmref{T3}}
We may assume that $X$ is locally given by $xy = tz$ at $p$.
Let $\pi: \wt{X}\to X$ be the blow-up of $X$ at $p$, i.e., $\wt{X}\subset
\A^3\times T\times
\P^3$ with coordinates $(x, y, z, t, x', y', z', t')$ be defined
by
$$
\left\{
\begin{array}{l}
x'y' = t'z',\\
\d{\frac{x}{x'}=\frac{y}{y'}=\frac{z}{z'}=\frac{t}{t'}}.
\end{array}
\right.
$$
Let $W = \wt{X}\cap \{t'\ne 0\}$. Let $u = x'/t'$ and $v = y'/t'$. Then 
$W\isom T\times \A_{uv}^2$ and $\pi: W\to X$ is given by
$$
\left\{
\begin{array}{lll}
x & = & tu,\\
y & = & tv,\\
z & = & tuv.
\end{array}
\right.
$$

Let $s_0, s_1, ..., s_n$ be a basis of $\sigma$ such that $s_1,s_2, ...,
s_{m^2}$ generates $(m-1)$-jets at $p$ on $X_0$. After choosing an
appropriate trivialization of $L$ and applying automorphisms of
$\sigma\times T$ induced by the action of $SL_{n+1}(\theta_t)$ on $\sigma$,
where $\theta_t$ is the ring of analytic power series in $t$ and
$SL_{n+1}(\theta_t)$ is the group of $(n+1)\times (n+1)$ matrices with
entries in $\theta_t$ and determinants 1, we have
\begin{align*}
\{s_1, s_2,..., s_{m^2}\} &=
\bigl\{x^iy^jz^k +
O\left(x^lz^{m-l}\big|_{l=0}^m, y^lz^{m-l}\big|_{l=0}^m
\right):\\
&\quad\quad i+j + k < m, ij = 0\bigr\}\\
s_i &= O\left(x^lz^{m-l}\big|_{l=0}^m,
y^lz^{m-l}\big|_{l=0}^m \right), \ {\mathrm{for}}\ m^2 < i \le n.
\end{align*}
Then the pullbacks of $s_1, s_2, ..., s_{m^2}, ..., s_n$ on $W$ are
\begin{align*}
\{\pi^*s_1, \pi^*s_2,..., \pi^*s_{m^2}\} &=
\bigl\{t^{\max(i,j)}u^iv^j+
O\left(t^mu^m, t^mv^m\right): \\
&\quad\quad i< m, j <m\bigr\}\\
\pi^*s_i &= O\left(t^mu^m, t^mv^m\right), \ {\mathrm{for}}\ m^2 < i \le n.
\end{align*}
Suppose the $2m$ branches of $C$ at $p$ are $\{\alpha_iy+z=0, x=0\}$ and
$\{\beta_ix+z = 0, y=0\}$ for $i=1,2,...,m$, where $\alpha_1\ne \alpha_2\ne
...\ne \alpha_m$ and $\beta_1\ne \beta_2\ne ...\ne \beta_m$. 
Without loss of generality, we may assume $s_0(x,y,z,t)$ satisfies
\begin{eqnarray*}
s_0(0, y, z, 0) &=& \prod_{i=1}^m (\alpha_iy+z) +
O\left(y^iz^j\big|_{i+j=m+1}\right),\\
s_0(x, 0, z, 0) &=& \prod_{i=1}^m (\beta_ix+z) +
O\left(x^iz^j\big|_{i+j=m+1}\right).
\end{eqnarray*}
It is not hard to see that there exists $c_1,c_2, ..., c_{m^2}\in t\theta_t$
such that
$$\pi^*s_0 - t^m\prod_{i=1}^m (\alpha_i + u)\prod_{i=1}^m (\beta_i + v)
= \sum_{i=1}^{m^2} c_i \pi^*s_i + O(t^{m+1}u^m, t^{m+1}v^m).$$
So we can apply an automorphism to $\sigma\times T$ such that
$$\pi^*s_0 = t^m\prod_{i=1}^m (\alpha_i + u)\prod_{i=1}^m (\beta_i + v) + 
O(t^{m+1}u^m, t^{m+1}v^m).$$
Let $B$ be a neighborhood of $s_0\in \sigma\isom \P^n$ and $Z\subset
B\times W$
be the family of curves over $B\times T$ whose fiber over $s\in B$ and
$t\in T$ is $\{s=0\}\cap W_t$, where $W_t$ is the fiber of $W\to T$ over
$t$. Then $Z$ is locally defined by
\begin{align*}
t^m\prod_{i=1}^m (\alpha_i + u)&\prod_{i=1}^m (\beta_i + v)
+\sum_{i=0}^{m-1}\sum_{j=0}^{m-1} t_{im+j+1}t^{\max(i,j)}u^iv^j\\
&+O\left(t_it^mu^m\big|_{i=1}^n, t_it^mv^m\big|_{i=1}^n,
t^{m+1}u^m, t^{m+1}v^m\right)=0
\end{align*}
where $(t_1, t_2,..., t_n)$ are coordinates of $B$ such that
$s=s_0+\sum_{i=1}^n t_is_i$ for $s\in B$.

Let $\psi: B'\times (T-\{0\})\to B\times (T-\{0\})$ be the base change
given by
\begin{eqnarray*}
t_{im+j+1} &=& t_{im+j+1}'t^{m-\max(i,j)}, \ {\mathrm{for}}\
0\le i < m, 0\le j <
m,\\
t_l &=& t_l', \ {\mathrm{for}}\ m^2 < l \le n,
\end{eqnarray*}
where $B'$ is an $n$-dimensional polydisk with coordinates $(t_1', t_2',
..., t_n')$. And let $Z'$ be the closure of $Z\times_{B\times T} \left(B'\times
(T-\{0\})\right)$ in $B'\times W$.
Then $Z'$ is locally defined by
\begin{eqnarray*}
\prod_{i=1}^m (\alpha_i + u)\prod_{i=1}^m (\beta_i + v)
&+&\sum_{i=0}^{m-1}\sum_{j=0}^{m-1} t_{im+j+1}'u^iv^j\\
&+&O\left(t_i'u^m\big|_{i>m^2}, t_i'v^m\big|_{i>m^2},
tu^m, tv^m\right)=0.
\end{eqnarray*}
It is obvious that the $m^2$ nodes on the general fiber of $Z'$ are the
deformations of $m^2$ nodes on the central fiber which is
$$\prod_{i=1}^m (\alpha_i + u)\prod_{i=1}^m (\beta_i + v) = 0.$$
The versal deformation space of a node $xy = 0$ is simply the family
$S\to\A^1$ where $S\subset \A^3$ is given by
$$xy+\tau = 0$$
and the map $S\to \A^1$ is the projection on $\tau$.
Let $\phi_{kl}$ be the natural map from $B'\times T$ to the versal
deformation space $\A^1$ of the node $u = -\alpha_k, v = -\beta_l$,
for $1\le k\le m$ and $1\le l\le m$.
Obviously, $\phi_{kl}^{-1}(\{0\})$ is the locus in $B'\times T$ where the
fiber of $Z'\to B'\times T$ has a node in the neighborhood of the point
$u = -\alpha_k, v = -\beta_l$. And it is easy to see that
$\phi_{kl}^{-1}(\{0\})$ is a smooth
hypersurface in $B'\times T$
with tangent plane $H_{kl}$ at the origin given by
$$\sum_{i=0}^{m-1}\sum_{j=0}^{m-1} t_{im+j+1}'(-\alpha_k)^i(-\beta_l)^j +
L_{kl}(t_{m^2+1}', ..., t_n')= 0,$$
where $L_{kl}(t_{m^2+1}', ..., t_n')$ is some linear combination of
$t_{m^2+1}', ..., t_n'$.
Obviously, these planes $H_{kl}$ ($1\le k\le m, 1\le l\le m$) intersect
transversely. Consequently,
$\cap \phi_{kl}^{-1}(\{0\})$ is a smooth subvariety of
$B'\times T$ with codimension $m^2$ given by
$$t_i' + O(t_j'\big|_{j>m^2}) = 0
\ {\mathrm{for}}\ 1\le i \le m^2.$$
Let
%$$Y' = \psi\left(\bigcap_{\genfrac{}{}{0pt}{}{1\le k\le m}{1\le l\le m}} \phi_{kl}^{-1}(\{0\})
%- \{(0, 0, ..., 0)\}\right)$$
$$Y' = \psi\left(\cap_{1\le k, l\le m} \phi_{kl}^{-1}(\{0\})
- \{(0, 0, ..., 0)\}\right)$$
and $\closure{Y'}$ be the closure of $Y'$ in $B\times T$.
Obviously, $Y'\subset Y$ and $\closure{Y'}$ is smooth with tangent plane
$t_1=t_2=...=t_{m^2}=0$ at the origin.
It remains to show that $Y' = Y$.

Suppose $\Y\subset X\to T$ be a one-parameter family of curves 
where $\Y_t\subset
X_t$ is cut out by an element of $\sigma$, $\Y_t$ has $m^2$ nodes in $U$
and $\Y_0 = C$. The $m^2$ nodes on each general fiber of $\Y\to T$
will give $m^2$ sections
after some base change $S\to T$ given by $t = s^n$. Suppose the $m^2$
sections are given by $(x_l(s), y_l(s), z_l(s))$ for $1\le l \le m^2$
accordingly. Then we claim all coordinates
$x_l(s), y_l(s)$ and $z_l(s)$ vanishes at $s=0$ with order
$n$. Apparently, the choice of the base change $S\to T$ is immaterial
here. All we want to say is that the $m^2$ nodes of $\Y_t$ approach point $p$
at the order of $t$. A more geometrical way to state this is as follows.

It is well known that a rational double point of a three-fold can be
resolved by blowing up at the point and blowing down along either ruling of
the exceptional quadric surface. Let $\pi_1: X_1\to X$ and 
$\pi_2: X_2\to X$ be the two corresponding resolutions of $X$ at $p$, where
the central fiber of $X_1\to T$ is a union of the blow-up $\wt{Q_1}$ of
$Q_1$ at $p$ and $Q_2$ by identifying the strict transform $\wt{E}_1$ of $E$ on
$\wt{Q_1}$ and $E$ on $Q_2$; symmetrically, 
the central fiber of $X_2\to T$ is a
union of $Q_1$ the blow-up $\wt{Q_2}$ of
$Q_1$ at $p$ and $Q_2$ by identifying $E$ on $Q_1$ and
the strict transform $\wt{E}_2$ of $E$ on $\wt{Q_2}$. 

Let $\Y_1=\pi_1^{-1}(\Y)$ be the
total transform of $\Y$ under $\pi_1: X_1\to X$. Then the central fiber
$\pi_1^{-1}(C)$ of
$\pi_1^{-1}(\Y)\to T$ is the union $\wt{C_1}\cup mF_1\cup C_2$, where
$\wt{C_1}$ is the proper transform of $C_1$ under the blow-up
$\wt{Q_1}\to Q_1$ and $F_1$ is the exceptional divisor of the
blow-up. Obviously, $\wt{C_1}$ meets $F_1$ at $m$ points $q_1, q_2, ...,
q_m$ corresponding to $m$ branchs of $C_1$ at $p$, $\{\alpha_iy+z = 0, x =
0\}$ for
$1\le i \le m$, and none of $q_1, q_2, ..., q_m$ lies on $\wt{E}_1$, i.e.,
none of $q_1, q_2, ..., q_m$ is the intersection $p' = F_1\cap C_2$. Then
we claim that the $m^2$ nodes on the general fibers of $\wt{\Y}\to T$
approach the $m$ points $q_1, q_2, ..., q_m$. More precisely, 
there are exactly $m$ nodes on the general fibers approaching
each point $q_i$, for $i=1,2,...,m$.
This actually implies our previous claim that $x_l(s)$ vanishes at $s=0$ with
order $n$ for $1\le l\le m^2$. 
The reason is quite staightforward. Let $X_1\subset X\times \P^1$ be given
by
$$\frac{x}{t}=\frac{z}{y}=\frac{\lambda}{\mu},$$
where $(\lambda,\mu)$ are the coordinates of $\P^1$. Under these
coordinates point $q_i$
is given by $x=y=z=0$ and $\lambda = -\alpha_i \mu$ for $i=1, 2, ..., m$. So
we actually have a more precise statement about the asymptotic behavior of
$x_l(s)$ which gives that 
for each $\alpha_i$, $1\le i\le m$, there are exactly $m$ sections
$(x_l(s),y_l(s), z_l(s))$ such that
$$x_l(s) = -\alpha_i s^n + \CO(s^{n+1}).$$
By symmetry, we have the same statement about the total transform
$\Y_2=\pi_2^{-1}(\Y)$ of $\Y$
under $\pi_2: X_2\to T$, which yields that for each
$\beta_j$, $1\le j\le m$, there are exactly $m$ sections
$(x_l(s), y_l(s), z_l(s))$ such that
$$y_l(s) = -\beta_j s^n + \CO(s^{n+1}).$$
In summary, we have
\begin{claim}\label{clm:1}
For each pair $(\alpha_i, \beta_j)$, $1\le i, j\le m$, there is
a unique section $(x_l(s), y_l(s), z_l(s))$ such that
$$
\left\{
\begin{array}{lll}
x_l(s) &=& -\alpha_i s^n + \CO(s^{n+1}),\cr
y_l(s) &=& -\beta_j s^n + \CO(s^{n+1}).
\end{array}
\right.
$$
\end{claim}

Let $\Upsilon$ be the nodal reduction of the family $\Y_1\to T$.
We have the diagram
\begin{equation}\label{E1}
\begin{array}{ccccccc}
\Upsilon & \mapright{} & \Y_1 & 
\subset & X_1 & \mapright{\pi_1} & X\\
\big\downarrow & & \big\downarrow & & & 
& \\
\Delta & \mapright{} & T & & & &\\
\end{array}
\end{equation}
Let $\eta: \Upsilon\to X_1$ be the morphism in \eqref{E1} and $\Gamma$ be
the central fiber of $\Upsilon\to \Delta$. We write $\Gamma$ as
$$\Gamma = \Gamma_0 \cup \Gamma_1 \cup \Gamma_2 \cup \Psi$$
where $\Gamma_0$ is the union of components on which $\eta$ is constant,
$\eta(\Gamma_1) = \wt{C_1}$, $\eta(\Gamma_2) = C_2$ and $\eta(\Psi)
= mF_1$.

Notice that the pull-back $\eta^*Q_2$ of the divisor $Q_2$ on $X$ can be
written as $\eta^*Q_2 = D_1+D_2$, where $\Supp(D_1) = \Gamma_2$ and
$\Supp(D_2) \subset \Gamma_0$. Since $F_1\cdot Q_2 = 1$, each irreducible
component of $\Psi$ must meet $\Gamma_2$. And since there is no 
connected component of $\Gamma_0$ disjoint from $\Gamma_1$,
$\Gamma_2$ and $\Psi$, we must have
$$p_a(\Gamma)\ge p_a(\Gamma_1)+p_a(\Gamma_2)-1,$$
where $p_a(\cdot)$ is the arithmetic genus of a curve.

Let $\wt{C_2}$ be the local normalization of $C_2$ in $U$. Obviously, we
have
$$p_a(\Gamma_1)\ge p_a(\wt{C_1}),\ p_a(\Gamma_2)\ge
p_a(\wt{C_2})$$
and
$$p_a(\wt{C_1}) = p_a(C_1) - \frac{m(m-1)}{2},\ p_a(\wt{C_2}) =
p_a(C_2) - \frac{m(m-1)}{2}.$$

Let $\omega$ be the dualizing sheaf of $Q_1\cup Q_2$. 
The adjunction formula produces
\begin{eqnarray*}
2p_a(C_1)-2 &=& \deg(K_{Q_1}\tensor L)\mid_{C_1}\\
&=& \deg(\omega|_{Q_1}(-E)\tensor L)\mid_{C_1}\\
&=& \deg(\omega|_{Q_1}\tensor L)\mid_{C_1} - E\cdot C_1\\
&=& \deg(\omega|_{Q_1}\tensor L)\mid_{C_1} - m,
\end{eqnarray*}
where $K_{Q_1}$ is the canonical line bundle of $Q_1$.
Similarly, we have
$$2p_a(C_2)-2 = \deg(\omega|_{Q_2}\tensor L)\mid_{C_2} - m.$$
And since
\begin{eqnarray*}
2p(C) - 2 &=& \deg(\omega\tensor L)\mid_C\\
&=& \deg(\omega|_{Q_1}\tensor L)\mid_{C_1} + 
\deg(\omega|_{Q_2}\tensor L)\mid_{C_2},
\end{eqnarray*}
we must have
\begin{equation}\label{E2}
p_a(\Gamma)\ge p_a(C) - m^2.
\end{equation}
Obviously, the general fiber $\Upsilon\to \Delta$ has arithmetic genus $p_a(C) -
m^2$. Hence the equality holds in \eqref{E2}, which happens only if
\begin{enumerate}
\item $\Gamma_1\isom\wt{C_1}$ and $\Gamma_2\isom\wt{C_2}$;
\item $\Gamma_1$ and $\Gamma_2$ are disjoint;
\item every connected component of $\Gamma_0\cup\Psi$ is a tree of
smooth rational curves and meets $\Gamma_1\cup\Gamma_2$ exactly once;
$\Gamma_0$ is contractible in the map $\eta: \Upsilon\to X_1$ and hence
$\Gamma_0 = \emptyset$.
\end{enumerate}
Furthermore, by \lemref{L:XY=T}, each component of $\Psi$ is
connected to some component of $\Gamma_2$ and hence $\Psi$ is disjoint from
$\Gamma_1$. Therefore, by \lemref{L:XmYn} (looking at the map $\Upsilon\to
\Y_1$), there are exactly $m$ nodes on the general fiber of $\Y_1\to T$
approaching each point $q_i$ for $i=1,2,...,m$. This finishes the proof of
\claimref{clm:1} and hence the proof of \thmref{T3}.

\end{proof}

\subsection{Standard deformation theory of planary curve singularities}

The basic deformation theorem on the first order
deformation of a curve on a smooth surface is Zariski's Theorem \cite{Z}.
Put in a form more suitable for our purpose, it says

\begin{thm}[Zariski]\label{t:zariski}
Let $W$ be any family of curves on a smooth surface $S$, $C$ be a
general member of $W$ and $C$ be reduced. By identifying the tangent space
$T_{[C]}W$ to $W$ at $C$ with a subseries of $H^0(N_{C/S})=H^0(\CO_C(C))$,
we have
\begin{enumerate}
\item $T_{[C]}W$ satisfies the adjoint condition of $C$, i.e.,
$T_{[C]}W\subset
H^0(\I\tensor\CO_C(C))$, where $\I$ is the adjoint ideal of $C$;
\item if $C$ has any singularities other than nodes, then actually
$T_{[C]}W\subset H^0(\J\tensor \CO_C(C))$ where $\J\subsetneqq\I$ is an ideal
strictly contained in $\I$.
\end{enumerate}
\end{thm}

As a direct application of Zariski's Theorem, \cite[Proposition 2.1,
p. 447]{H} gives the upper bound of the dimension of
a family of reduced curves with fixed geometric genus on a rational surface.
This can be generalized in various ways. For
example, we can further impose some tangency conditions on the family of
curves. Specifically, we have

\begin{thm}\label{T:HARRIS}
Let $S$ be a smooth rational surface, $D$ a divisor class on $S$,
$W\subset |D|$ a family of reduced curves of geometric genus $g$ and
$C\subset S$ be a reduced curve. Suppose that $C$ meets a general member
$E$ of $W$ at $r$ points $p_1, p_2, ..., p_r$, 
which are smooth on both $C$ and $E$,
with multiplicities $m_1, m_2, ..., m_r$, respectively. 
And suppose that for a general member $E\in W$, the restriction of the
divisor $-(K_S + C)|_E + \sum_{i=1}^r p_i$ to each irreducible component of
$E$ has degree at least 2, where $K_S$ is the canonical
divisor of $S$.

Let $\delta =
-D\cdot K_S - \sum_{i=1}^r (m_i-1)$. We have
\begin{enumerate}
\item $\dim W \le \max(\delta + g - 1, 0)$;
\item if $\dim W = \delta + g - 1 > 0$, the general member of $W$ is nodal;
\item if $\dim W = \delta + g - 1 > 1$ and $F\subset S$ is a reduced
curve intersecting $C$ properly, the general member of $W$ meets $F$
transversely.
\end{enumerate}
\end{thm}

\begin{rem}
Different versions of this theorem have already appeared in \cite{H},
\cite{CH1} and \cite{CH2}, but all with slightly different hypotheses, for
example, with a single tangency condition to a line in \cite[Lemma 2.4,
p. 450]{H}, dealing with $S = \F_n$ in \cite[Proposition 2.1]{CH1} and
dealing with plane curves in \cite[Proposition 2.1]{CH2}. But all
techniques necessary to prove this theorem can be found in those places.

Basically, one notices that the tangent space $T_{[E]} W$ (assuming $E$ is
irreducible) is contained in $H^0(\J\tensor \CO_E(E)\tensor
\CO_E(-\sum_{i=1}^r (m_i-1) p_i))$ by Zariski's theorem, where 
$C\cdot E = \sum_{i=1}^r m_i p_i$ and
$\CO_E(-\sum_{i=1}^r (m_i-1) p_i)$ accounts for the tangency conditions.
Again by Zariski's theorem, $\J\tensor \CO_E(E)\tensor
\CO_E(-\sum_{i=1}^r (m_i-1) p_i)\subsetneqq \I\tensor \CO_E(E)\tensor
\CO_E(-\sum_{i=1}^r (m_i-1) p_i)$.

Let $v: E^v\to E$ be the normalization of $E$. We have 
\[
\begin{split}
& \quad v^* (\I\tensor
\CO_E(E)\tensor \CO_E(-\sum_{i=1}^r (m_i-1) p_i)) \\
&= K_{E^v}\tensor
v^* K_S^{-1} \tensor \CO_{E^v}(-\sum_{i=1}^r (m_i-1) p_i)\\
&= K_{E^v}\tensor \CO_{E^v} (-K_S - C) \tensor \CO_{E^v} (\sum_{i=1}^r p_i).
\end{split}
\]
Our numerical conditions on $S$, $C$ and $E$ guarantee that the complete
linear series $|K_{E^v}\tensor
\CO_{E^v} (-K_S - C) \tensor \CO_{E^v} (\sum_{i=1}^r p_i)|$ is base point
free on $E^v$. And thus the argument in \cite{H} applies here.

Also notice the dimension requirement on $W$ in order to conclude that the
general member of $W$ is nodal or meets a fixed curve transversely. This
point is not stressed in the places mentioned above since the expected
dimension $\delta + g - 1$ is big enough there. However, in our
application, this is essential since we are dealing with cases $\delta + g
- 1\le 2$.
\end{rem}

\subsection{Review of general deformation theory}

It will come up in our attempt to
the proof of \conjref{T:NOD} that we need to study the
deformation of a trigonal K3 surface in the projective space. Hence we will
give a review of some simple aspects of deformation theory. For our
purpose, we will only concern ourselves with embedded deformations.

Let $X$ be a closed subscheme of $Y$. An
embedded first order deformation of $X$ in $Y$ is a scheme $W\subset
Y\times D$ which is flat over $D$ with central fiber $X$ where $D=\Spec
\BC[t]/(t^2)$. Such $W$'s are classified by $H^0(N_{X/Y})$ as follows.

Let $\I_X$ be the ideal sheaf of $X$. A global section $s$ of the normal
sheaf $N_{X/Y})$ gives a sheaf morphism $\I_X/\I_X^2\to \CO_X$. Let $\Spec
A$ be an affine open set of $Y$ and $I_X = \Gamma(\Spec A, \I_X)$. Then
the embedded first order deformation $W$ corresponding to $s$ is locally
defined by the ideal generating by $f+tg$ where $f\in I_X$, $g\in A$ and
$s(f) = g$ when $f$ and $g$ are restricted to $I_X/I_X^2$ and $A/I_X$,
respectively.

Of course, we are really interested in the deformations of $X$ in $Y$ over
disk $T$. Studying the first order deformations is the first step to
classify deformations over disk $T$. Then it raises a natural question that
when a first order deformation can be ``lifted'' to a deformation over
$T$. By a ``lift'' we mean a scheme $V\subset Y\times T$
flat over $T$ such that $W\isom V\times_T D$.
Suppose $Y$ is smooth and $X$ is a locally complete intersection, which is
satisfied in our case.
The obstruction to lift a first order deformation turns out to be
$H^1(N_{X/Y})$. If $H^1(N_{X/Y}) = 0$,
every first order deformation can be lifted. Moreover, if $H^1(N_{X/Y})$
vanishes, any embedded deformation of $X$ over $\Spec \BC[t]/(t^l)$, i.e.,
a scheme $W_l \subset Y\times \Spec \BC[t]/(t^l)$ flat over $\Spec
\BC[t]/(t^l)$ with central fiber $X$, can be lifted to a deformation
over $T$, i.e., a scheme $V\subset Y\times T$ flat over $T$
satisfying $W_l\isom V\times_T \Spec \BC[t]/(t^l)$.

\section{Existence of Rational Curves on a General K3 Surface}\label{S:EXIST}

\subsection{Degeneration of K3 surfaces}
A general K3 surface can be degenerated to a union of two rational scrolls.
For example, for a quartic surface $S$ in $\P^3$, we may simply take the
degenerating family as the pencil connecting $S$ with the union of two smooth 
quadric surfaces $Q_1\cup Q_2$ in general position.

In general, it was shown in \cite{CLM} that the union of two
rational normal scrolls (each of degree $n-1$ in $\P^n$) meeting
transversally along a smooth anticanonical elliptic curve lies on the
boundary of the component of Hilbert scheme consisting of primitive
K3 surfaces in $\P^n$, i.e., K3 surfaces in $\P^n$ on which $\CO(1)$ is
non-divisible.
Here is a sketch of their proof.

Let $R = R_1\cup R_2$ denote the union of two rational normal scrolls of
degree $n-1$ in $\P^n$ and
$E = R_1\cap R_2$ be the smooth elliptic curve cut out by $R_1$ and $R_2$. 
Let $T^1 = {\mathcal
E}xt^1(\Omega_R, \CO_R)$ and $N_R$ be the normal bundle of $R$ in $\P^n$. It
is shown in \cite[Sec. 2.2]{CLM} that $H^1(N_R)$ vanishes and
$H^0(N_R)$ surjects onto $H^0(T^1)$. By a
standard deformation theorem \cite{F}, the embedded deformations of $R$ in
$\P^n$ smooth the double curve of $R$ and hence deform $R$ to a $K3$
surface in $\P^n$.

Furthermore, since $T^1$ is a coherent sheaf supported on $E$
whose restriction on $E$ is the line bundle
$N_{E/R_1}\tensor N_{E/R_2}$, it is easy to see that a general one-parameter
family of K3 surfaces with central fiber $R$ has exactly 
16 distinct ordinary double points in general position on $E$.

Let $X\subset\P^n\times T$
be a general one-parameter family of K3 surfaces over disk $T$
whose central fiber $X_0 = Q_1\cup Q_2\subset\P^n$
is a union of rational normal scrolls of degree $n-1$.
And $Q_1$ and $Q_2$ meet transversely along
a smooth anticanonical elliptic curve $E$. 
Let $p_1, p_2, ..., p_{16}$ be the sixteen rational double points 
of $X$ and $l=\lfloor n/2 \rfloor$.

Notice that since $Q_i$ contains a smooth elliptic curve, $Q_i$ must be
$\P^1\times \P^1$, $\F_1$ or $\F_2$, where $\F_m = \P(\CO\oplus\CO(-m))$ is
the rule surface.

For $n$ odd, $Q_i$ can either be $\P^1\times\P^1$ or $\F_2$. We choose
$Q_i$ to be $\P^1\times \P^1$ embedded into $\P^n$ by the divisor
$H_1+lH_2$ for $i=1,2$, where
$H_1 = \P^1\times\{pt\}$ and $H_2 = \{pt\}\times\P^1$. 

For $n$ even, we must have
$Q_i\isom\F_1$ embedded into $\P^n$ by the divisor 
$C + lF$ for $i=1,2$,
where $C$ is the divisor associated to the line bundle
$\CO_{\P(\CO\oplus\CO(-1))}(1)$ and $F$ is a fiber of
the projection $\P(\CO\oplus\CO(-1))\to \P^1$.

To prove the existence theorem \ref{T:EXIST},
it suffices to locate a limiting rational curve in the linear series
$|\CO(d)|$ on the central fiber
$Q_1\cup Q_2$ of the family $X$ constructed above. 

\subsection{The curve we are looking for}
We are looking for a curve $C_{1}^{1}\cup C_{2}^{1}\cup ...\cup
C_{d-1}^{1}\cup C_{d}^{1}\cup C_{1}^{2}\cup C_{2}^{2}\cup ...\cup
C_{d-1}^{2}\cup C_{d}^{2}$ on the central fiber $X_0=Q_1\cup Q_2$ of the
degenerating family $X$ constucted at the beginning of this section, where
\begin{enumerate}
\item $C_j^{i}\subset Q_i$ for $i=1,2$ and $1\le j\le d$;
\item $C_j^{i}\in |H_1|$ for $1\le j \le d-1$
and $C_{d}^{i}\in |H_1+dlH_2|$
if $n$ is odd;
$C_j^{i}\in |C+F|$ for $1\le j \le d-1$ and $C_{d}^{i}\in |C+(dl-d+1)F|$
if $n$ is even, for $i=1,2$;
\item these curves $C_i^j$ ($1\le i\le 2, 1\le j\le d$) are determined by the
following relations (let $q_0=p_1$)
$$
\begin{array}{ll}
C_{j}^{1}\cap E = q_{2j-2}+q_{2j-1},&
C_{j}^{2}\cap E = q_{2j-1}+q_{2j},\ {\mathrm{for}}\ j <  d,\\
C_{d}^{1}\cap E = q_{2d-2}+ (2dl+1)r, & C_{d}^{2}\cap E = q_0+ (2dl+1)r,
\end{array}
$$
if $n$ is odd and
$$
\begin{array}{ll}
C_{j}^{1}\cap E = q_{2j-2}+2q_{2j-1}, &
C_{j}^{2}\cap E = 2q_{2j-1}+q_{2j},\ {\mathrm{for}}\ j < d,\\
C_{d}^{1}\cap E = q_{2d-2}+ (2dl-2d+2)r,
& C_{d}^{2}\cap E = q_0+ (2dl-2d+2)r,
\end{array}
$$
if $n$ is even, where $q_1,q_2,...,q_{2d-2}$ and $r$ are points on $E$.
\end{enumerate}
We see that the points $q_1,q_2,...,q_{2d-2}$ and $r$ are uniquely determined
by these relations. Since $q_0=p_1,p_2,...,p_{16}$ are in general position
on $E$, and the embeddings $i_1: E\embed Q_1$ and $i_2: E\embed Q_2$ are
general under the condition ${i_1}^*\CO_{Q_1}(H_1+lH_2) = {i_2}^*
\CO_{Q_2}(H_1+lH_2)$ if $n$ is odd and ${i_1}^*\CO_{Q_1}(C+lF) = {i_2}^*
\CO_{Q_2}(C+lF)$ if $n$ is even, we may assume $q_1,q_2,...,q_{2d-2}$ and
$r$ are different from each other and not among $p_1,p_2,...,p_{16}$. We
may further assume the curves $C_{j}^{i}$ intersect each other transversely
and no three of them meet at a point.

A straight calculation shows the complete linear series of 
$\CO_{\P^n}(d)\big|_{Q_1\cup Q_2}$ consists of exactly the curves $C_1\cup C_2$ 
where
\begin{enumerate}
\item $C_1\subset Q_1$ and $C_2\subset Q_2$;
\item $C_i\in |dH_1+dlH_2|$ if
$n$ is odd and $C_i\in |dC + dlF|$ if $n$ is even, for $i=1,2$;
\item $C_1\cap E = C_2\cap E$.
\end{enumerate}
Hence $\cup C_{j}^{i}$ is cut out by a hypersurface of degree $d$ in $\P^n$.

This curve may look strange at first. But geometrically it is
quite clear how $\cup C_{j}^{i}$ can be deformed
to a rational curve on the general fiber. For example, 
in the case that $n$ is odd,
point $r$ deforms to $2dl$ nodes by \thmref{T2},
point $p_1$ deforms to a node by \thmref{T3}. 
And the intersections, $C_{j}^{i}\cap C_{d}^{i}$ 
($1\le i\le 2, 1\le j\le d-1$), can be deformed equisingularly.

\subsection{Completion of the proof of the existence theorem}

We will only finish the proof for the case $n$ is odd,
since the same argument applies to the case $n$ is even almost without
change.

Let $U^{d,\delta}(S)$ be the subset of $|\CO_S(d)|$ consisting of
irreducible nodal curves with $\delta$ nodes on a K3 surface $S\subset\P^n$.
Let
$$
\begin{array}{ll}
Y^{d, \delta}(S) = & \{ (C, x_1, x_2, ..., x_\delta, t): C\in
U^{d,\delta}(S),\\
& x_1, x_2, ..., x_\delta\in\P^n \ {\mathrm{are\ distinct\ nodes\ of}}\ C, 
t\ne 0\in T\}.
\end{array}
$$
And let $Y^{d, \delta}(X)$ be the fiberation over $T$ whose fibers are
$Y^{d,\delta}(X_t)$, where $Y^{d, \delta}(X_0)$ is the flat limit of
$Y^{d, \delta}(X_t)$.

Let $\pi_1: Y^{d,\delta}(X)\to\P^{d^2(n-1)+1}$, 
$\pi_k: Y^{d,\delta}(X)\to\P^n$ for $2\le k\le \delta+1$ and 
$\pi_{1,\delta+2}: Y^{d,\delta}(X) \to \P^{d^2(n-1)+1}\times T$ be
projections from $Y^{d,\delta}(X)$ as a subscheme of 
$\P^{d^2(n-1)+1}\times(\P^n)^\delta\times T$.

Let $C=\cup C_{j}^{i}$ and 
$x_1, x_2,..., x_{2(d-1)dl}$ be the intersections among curves
$C_{j}^{i}$ not lying on $E$. 
Obviously $(C, x_j)\in Y^{d, 1}(X_0)$
for all $1\le j\le 2(d-1)dl$. By \thmref{T2}, $(C, r, ..., r)\in
Y^{d, 2dl}(X_0)$ and by \thmref{T3}, $(C, p_1)\in Y^{d, 1}(X_0)$. 

Choose analytic
neighborhoods $O_j$ of $(C, x_j, 0)$ in $Y^{d,1}(X)$ for $1\le
j\le 2(d-1)dl$ 
such that $\pi_2(O_1), \pi_2(O_2), ..., \pi_2(O_{2(d-1)dl})$ are disjoint
from each other and the curve $E$. Obviously, $\codim\pi_{1,3}(P)\le 1$ in
$\P^{d^2(n-1)+1}\times T$ for $1\le j\le 2(d-1)dl$. By \thmref{T3},
there is also a neighborhood $P$ of $(C, p_1, 0)$ in $Y^{d,1}(X)$ such that 
the central fiber of $\pi_{1,3}(P)$ over $T$ consists of curves passing
through point $p_1$ and hence $\codim\pi_{1,3}(P)\le 1$ in
$\P^{d^2(n-1)+1}\times T$. 
We may make $\pi_2(P)$
disjoint from each of $\pi_2(O_1), \pi_2(O_2), ..., \pi_2(O_{2(d-1)dl})$ and
point $r$. 
Similarly, by \thmref{T2},
We can also choose a
neighborhood $Q$ of $(C, r, ..., r, 0)$ in $Y^{d, 2dl}(X)$ such that
the central fiber of $\pi_{1, 2dl+2}(Q)$ over $T$
consists of curves meeting $E$ at
$r$ with multiplicity $2dl+1$, $\codim \pi_{1,2dl+2}(Q)\le 2dl$ in
$\P^{d^2(n-1)+1}\times T$ and
$\pi_2(Q), \pi_3(Q), ..., \pi_{2dl+1}(Q)$ are disjoint from each of
$\pi_2(O_1), \pi_2(O_2), ..., \pi_2(O_{2(d-1)dl}), \pi_2(P)$.

Let 
$W = \cap \pi_{1,3}(O_j)\cap \pi_{1,3}(P)\cap \pi_{1,2dl+2}(Q)$. It is easy
to see the central fiber of $W$ over $T$ consists of curve $C$ and the
general fiber of $W$ over $T$ consists of curves with at least $2d^2l+1$
nodes. Since $\dim W \ge 1$, the general fiber of $W$ over $T$ is nonempty.
Consequently,
there exists a family
of curves $C_t$ over $T$ such that $C_0 = C$ and $C_t\in |\CO_{X_t}(d)|$ has
at least $2d^2l+1$ nodes. Besides,
it is not hard to see the general fiber $C_t$ is irreducible. Otherwise, if
$C_t$ contains a curve in $|\CO_{X_t}(d')|$ for some $d'<d$, then $C_0=C$
must constain a curve in $|d'(H_1+lH_2)|$ on $Q_i$. This contradicts
the choice of $C=\cup C_j^i$. Hence $C_t$ must be an irreducible nodal
rational curve. This finishes the proof of \thmref{T:EXIST}.

A generalization of \thmref{T:EXIST} can be made on the curves
on a general K3 surface with any given geometric genus.

\begin{thm}\label{T:EXISTA}
$U^{d,\delta}(S)\ne \emptyset$ 
and $\codim U^{d,\delta}(S)
= \delta$ on a general K3 surface $S$ in $\P^n$ 
for each $\delta\le d^2(n-1)+1$.
\end{thm}

Given \thmref{T:EXIST}, it suffices to bound $\dim U_{d,g}(S)$ from above
in order to prove \thmref{T:EXISTA}.

\begin{lem}\label{L1}
Let $U_{d,g}(S)$ be the subset of
$|\CO_S(d)|$ consisting of reduced and irreducible curves of geometric genus
$g$ on a K3 surface $S\subset\P^n$. Then
$\dim U_{d,g}(S) \le g$. If $g > 0$ and $W$ is an irreducible component of
$U_{d,g}$ of dimension $g$, then the general member of $W$ is a nodal curve.
\end{lem}

This is an easy application of Zariski's Theorem.

\begin{proof}{Proof of \lemref{L1}}
Let $W$ be an irreducible component of $U_{d,g}(S)$, $C$ be a general
member of $W$ and $v: C^v\to C$ be the normalization of $C$.
Since $\I$ imposes independent conditions on $|\omega_C|=|K_S\tensor
\CO_C(C)| = |\CO_C(C)|$, by Zariski's Theorem we have
\begin{eqnarray*}
\dim T_{[C]}W & \le & h^0(\I\tensor\CO_C(C))\\
& = & h^0(\CO_C(C)) - (p_a(C) - g) = g,
\end{eqnarray*}
where $p_a(C)$ is the arithmetic genus of $C$.
Therefore, $\dim W \le g$.

Since $v^*H^0(\I\tensor\CO_C(C))$ cuts out the complete series 
$|K_{C^v}\tensor v^* K_S^{-1}| = |K_{C^v}|$ on $C^v$, 
which is base point free for $g>0$, $H^0(\I\tensor\CO_C(C))$ must be base
point free on $C$ for $g>0$.
Hence if $\dim W = g > 0$ and $C$ has singularities other
than nodes, we must have
$$\dim W \le H^0(\J\tensor\CO_C(C)) <
H^0(\I\tensor\CO_C(C)) = g.$$
Contradiction.
\end{proof}

\section{Degeneration To Trigonal K3 Surfaces}\label{S:TRIG1}

We will spend the rest of the paper showing the progress we have made
towards \conjref{T:NOD}.
Though the degeneration of a K3 surface to a union of rational
normal scrolls helps to establish the existence theorem, it fails here
due to the presence of nonreduced limiting rational curves for
$n\ge 5$. An alternative degeneration will be
introduced. Basically, we will do the
degeneration in two steps. First we degenerate a general K3 surface to a
trigonal K3 surface. Then we further degenerate a trigonal K3 surface to a
union of rational surfaces. The rest of this paper will concentrate on the
first step of this degeneration.

A trigonal K3 surface in $\P^n$ is a K3 surface containing a pencil of
elliptic curve of degree 3, namely, a K3 surface with Picard lattice
congruent to
$\bigl(\begin{smallmatrix}
2n-2 & 3\\
3 & 0
\end{smallmatrix}\bigr).$
The transcendental theory of
K3 surfaces shows that the moduli space of trigonal K3 surfaces consists of
countably many irreducible components of dimension 18. We need three of them

\begin{description}
\item[{\bf{TK1}}] surfaces in $\P^2\times \P^1$ of type $(3, 2)$
embedded into $\P^n$ by the line bundle $\CO(1, k)$ for $n=3k+2$ ($k>0$);
\item[{\bf{TK2}}] complete intersections of $(3,1)$ and $(1,1)$ hypersurfaces in
$\P^3\times\P^1$ 
embedded into $\P^n$ by the line bundle $\CO(1,k)$ for $n=3k+3$ ($k>0$);
\item[{\bf{TK3}}] complete intersections of $(3,0)$, $(1,1)$ and $(1,1)$ 
hypersurfaces in
$\P^4\times\P^1$ 
embedded into $\P^n$ by the line bundle $\CO(1,k)$ for $n=3k+4$ ($k>0$).
\end{description}

Alternatively, we can think of these surfaces as the anticanonical surfaces
of projective bundles $\P E$ over $\P^1$ embedded into $\P^n$ by
$\CO(C+kF)$, where $E = \CO\oplus\CO\oplus\CO$,
$\CO\oplus\CO\oplus\CO(1)$ or $\CO\oplus\CO(1)\oplus\CO(1)$ corresponding
to {\bf{TK1}}, {\bf{TK2}} or {\bf{TK3}}, respectively, and
$C$ and $F$ are the divisors on $\P E$ corresponding to the line bundles
$\CO_{\P E}(1)$ and $\pi^*\CO_{\P^1}(1)$ ($\pi: \P E\to \P^1$).

\begin{prop}\label{P1}
Let $S$ be a K3 surface given in {\bf{TK1}}, {\bf{TK2}} or 
{\bf{TK3}} and $N_S$ be the
normal bundle of $S$ in $\P^n$. Then $\dim H^0(N_S) = n^2+2n+19$ and 
$H^1(N_S) = 0$.
\end{prop}

\begin{proof}
Let $N_{\P E}$ be the normal
bundle of $\P E$ in $\P^n$. We have the exact sequence
\begin{equation}
0\mapright{} N_{S/\P E}\mapright{} N_S \mapright{}
N_{\P E}|_S\mapright{} 0.
\end{equation}
Obviously, $H^0(N_{S/\P E}) = 29$ and $H^1(N_{S/\P E}) = 0$. Hence
it suffices to show that $\dim H^0(N_{\P E}|_S) = n^2+2n-10$ and
$H^1(N_{\P E}|_S) = 0$.
Fixing a section $\P^1\to\P E$, we have 
\begin{equation}\label{e:p1-1}
T_{\P E} = T_{\P E/\P^1}\oplus
\pi^*T_{\P^1},
\end{equation}
where $\pi: \P E\to\P^1$ is the projection.
And we have the Euler sequences
\begin{equation}\label{e:p1-2}
0 \mapright{} \CO \mapright{} \CO(1) \oplus \CO(1) \mapright{} T_{\P^1} \mapright{} 0
\end{equation}
on $\P^1$,
\begin{equation}\label{e:p1-3}
0\mapright{} \CO \mapright{} \pi^* E^\vee \tensor \CO(C) \mapright{} T_{\P E
/\P^1} \mapright{} 0
\end{equation}
over $\P E\to \P^1$ and
\begin{equation}\label{e:p1-4}
0\mapright{} \CO \mapright{} \CO(1)^{\oplus (n+1)} \mapright{} T_{\P^n}
\mapright{} 0
\end{equation}
on $\P^n$.
Combining \eqref{e:p1-1}, \eqref{e:p1-2} and \eqref{e:p1-3} and restricting
them to $S$, we have
\begin{equation}\label{e:p1-5}
0 \mapright{} \CO_S\oplus\CO_S \mapright{}
\pi^*E^\vee\tensor\CO_S(C)\oplus \CO_S(F)^{\oplus 2} \mapright{} 
T_{\P E}|_S \mapright{} 0.
\end{equation}
And restricting \eqref{e:p1-4} to $S$, we have
\begin{equation}\label{e:p1-6}
0 \mapright{} \CO_S \mapright{} \CO_S(C+kF)^{\oplus (n+1)} \mapright{}
T_{\P^n}|_S \mapright{} 0.
\end{equation}
Since $\Ext^1(\pi^*E^\vee\tensor\CO_S(C)\oplus \CO_S(F)^{\oplus 2}, \CO_S)
= H^1(\pi^* E \tensor\CO_S(-C)\oplus \CO_S(-F)^{\oplus 2}) = 0$, there
exists $v\in \Hom(\pi^*E^\vee\tensor\CO_S(C)\oplus \CO_S(F)^{\oplus 2},
\CO_S(C+kF)^{\oplus (n+1)})$ such that the diagram
\begin{equation}
\begin{array}{ccccccccc}
& & 0 & & 0\\
& & \downarrow & & \downarrow\\
0 & \to & \CO_S & \to & \CO_S & \to & 0\\
 & & \downarrow & & \downarrow & & \downarrow\\
0 & \to & \CO_S\oplus\CO_S & \to &
\pi^*E^\vee\tensor\CO_S(C)\oplus \CO_S(F)^{\oplus 2} & \to & 
T_{\P E}|_S & \to & 0 \\
& & \downarrow & & \smapdown{v} & & \downarrow & & \\
0 & \to & \CO_S & \to & \CO_S(C+kF)^{\oplus (n+1)} & \to
& T_{\P^n}|_S & \to & 0\\
& & \downarrow & & \downarrow & & \downarrow & &\\
& &  0 & \to & N_{\P E}|_S & \to &
N_{\P E}|_S &\to & 0\\
& & & & \downarrow & & \downarrow\\
& & & & 0 & & 0\\
\end{array}
\end{equation}
is commutative and exact in each column and row.
We are interested in the middle column
\begin{equation}
\begin{split}
0\mapright{} \CO_S & \mapright{} \pi^*E^\vee\tensor
\CO_S(C)\oplus \CO_S(F)^{\oplus 2}\\
& \mapright{v} \CO_S(C+kF)^{\oplus (n+1)} \mapright{} N_{\P E}|_S
\mapright{} 0.
\end{split}
\end{equation}
Since $H^i(\pi^*E^\vee\tensor\CO_S(C)\oplus \CO_S(F)^{\oplus 2}) = H^i(\CO_S(C+
kF)^{\oplus (n+1)}) = 0$ for $i>0$, $H^1(N_{\P E}|_S) = H^3(\CO_S)
= 0$, $H^2(N_{\P E}|_S) = 0$ and
\begin{eqnarray*}
h^0(N_{\P E}|_S) & = & \chi(\CO_S) + h^0(\CO_S(C+kF)^{\oplus (n+1)}) \\
&& - h^0(\pi^*E^\vee\tensor\CO_S(C)\oplus \CO_S(F)^{\oplus 2})\\
& = & n^2+2n-10.
\end{eqnarray*}
\end{proof}

It follows \propref{P1} that a K3 surface $S$ in {\bf{TK1}}, {\bf{TK2}} or 
{\bf{TK3}} represents a smooth point of the Hilbert
scheme of K3 surfaces in $\P^n$ and it lies on an irreducible component
of dimension $n^2+2n+19$. And since the restriction of $\CO_{\P^n}(1)$ to
$S$ is indivisible, by the transcendental theory of K3 surfaces
$S$ lies on the component of the Hilbert
scheme consisting of primitive K3 surfaces in $\P^n$.

Again we take a one-parameter family of general K3 surfaces whose
central fiber $S$ is a trigonal K3 surface given in {\bf{TK1}}, {\bf{TK2}} or 
{\bf{TK3}} and we ask which curves on $S$ are limiting rational curves. Let
$\Gamma$ be a limiting rational curve in the form
\begin{equation}\label{E3}
\Gamma = \Gamma_0\cup m_1 \Gamma_1 \cup m_2 \Gamma_2\cup ...\cup 
m_\alpha\Gamma_\alpha
\end{equation}
where $\Gamma_i$ ($0\le i\le \alpha$) are irreducible components of
$C$ with multiplicities $m_i$ (let $m_0=1$), $\Gamma_0\in |\CO_S(C+lF)|$ and
$\Gamma_i \in |\CO_S(F)|$ for $i>0$ where $l+\sum_{i>0} m_i =
k$. Obviously, being a limit of rational curves, $\Gamma_i$ must be
rational. We will not go into the study of rational curves on trigonal K3
surfaces in this paper. Instead, we will assume the following

\begin{conj}\label{T:NODA}
Let $S$ be a surface given in {\bf{TK1}}, {\bf{TK2}} or 
{\bf{TK3}}. Then every irreducible rational curve in $|F|$ is nodal and
it intersects transversely with any irreducible rational curve in
$|C+lF|$.
\end{conj}

In this section we will show that 

\begin{thm}\label{T:NODB}
Suppose that \conjref{T:NODA} is true.
Let $W\subset \P^n\times T$ be a family of K3 surfaces over disk $T$ whose
central fiber $S$ is a surface given in {\bf{TK1}}, {\bf{TK2}} or 
{\bf{TK3}}.
Let $\Upsilon\subset W$ be a family of rational curves cut out by
$H^0(\CO_{\P^n}(1))$ with central fiber $\Gamma$ in the form \eqref{E3}. Then
$m_1 = m_2 = ... = m_\alpha = 1$, i.e., $\Gamma$ is reduced for $W$
general.
\end{thm}

If both \conjref{T:NODA} and \thmref{T:NODB} hold, 
\conjref{T:NOD} is an immediate consequence of the following statement

\begin{conj}\label{T:NODC}
Let $S$ be a surface given in {\bf{TK1}}, {\bf{TK2}} or 
{\bf{TK3}}. Then every irreducible rational curve in $|C+lF|$
is nodal.
\end{conj}

Since the arguments for the three cases {\bf{TK1}}, {\bf{TK2}} and
{\bf{TK3}} are essentially similar to each other,
we will only deal with {\bf{TK1}} here. 

Also note that this degeneration only works for $n\ge 5$. While for
$n<5$ we can work out \conjref{T:NOD} in a straightforward way as
follows.

For $n=3$, let $W\subset  |\CO_{\P^3}(1)|\times
|\CO_{\P^3}(4)|$ be the incidence correspondence
$(H, S)$ such that
$H\cap S$ is an irreducible rational curve.
Projecting $W$ to $\P^3$, we see that a fiber of $W$
over $H\in |\CO_{\P^3}(1)|$ can be identified with
$V_{4, 0}\times H^0(\CO_{\P^3}(3))$, where $V_{d,0}$ is
the Severi variety of degree $d$ irreducible rational curves on $\P^2$,
which is irreducible. Hence $W$ is irreducible and
we actually have \conjref{conj:1} for $n=3$ and
$d=1$.

For $n=4$, it is well-known that
every K3 surface in $P^4$ is a complete intersection of a
quadric and a cubic. Let $W\subset |\CO_{\P^4}(1)|\times |\CO_{\P^4}(2)|
\times |\CO_{\P^4}(3)|$ be the
incidence correspondence $(H, Q, C)$ such that 
$H\cap Q\cap C$ is an irreducible rational curve.
Projecting $W$ to $\P^4$ as above, we see that
a fiber of $W$ over $H\in |\CO_{\P^4}(1)|$ can be identified with
$V\times H^0(\CO_{\P^4}(1))\times
H^0(\CO_{\P^4}(2))\times |\CO_{\P^3}(2)|\times H^0(\CO_{\P^3}(1))$,
where $V$ is
the variety parameterizing irreducible
rational curves of type $(3,3)$ on $\P^1\times\P^1$, which is irreducible.
Hence $W$ is irreducible and we actually have \conjref{conj:1} for $n=4$
and $d=1$.

\subsection{Deformation of a Trigonal K3 surface}

Let $S$ be a $(3,2)$ surface in $\P^2\times\P^1$ and the embedding of
$\P^2\times\P^1$ into $\P^{3k+2}$ be given by
$$Z_{ij} = X_jY_0^{k-i}Y_1^i$$
where $(X_0, X_1, X_2)\times (Y_0, Y_1)$ and $(Z_{ij})$ ($0\le i\le 2, 0\le
j\le k$) are the projective coordinates of
$\P^2\times\P^1$ and $\P^{3k+2}$, respectively.
Also let 
$x_1=X_1/X_0$, $x_2=X_2/X_0$ and $y_1=Y_1/Y_0$ be the affine coordinates of
$\P^2\times\P^1$ over the open set $X_0Y_0\ne 0$. And correspondingly let
$z_{ij} = Z_{ij}/Z_{00}$ be the affine coordinates of $\P^{3k+2}$ over the
open set $Z_{00}\ne 0$. Then the embedding of $\P^2\times\P^1$ into
$\P^{3k+2}$ is locally given by
$$z_{ij} = x_jy_1^i.$$
Let $S$ be defined by
\begin{equation}\label{E32}
f(x_1, x_2, y_1) = q + x_1h_1 + x_2h_2 = 0
\end{equation}
where $q\in \BC[y_1]$ is a quadratic polynomial in $y_1$.
Without loss of generality, let us assume that points $(X_0=X_1=Y_1=0)$ and
$(X_0=X_2=Y_1=0)$ do not lie on $S$.
Obviously,
$\P^2\times \P^1$ is a local complete intersection in $\P^{3k+2}$ and so is
$S$. Take $z_{ij} - x_j y_1^i$ ($i\ne 0$ and $(i,j)\ne(1,0)$)
as the defining polynomials of
$\P^2\times \P^1$ in $\P^{3k+2}$ (here we identify $x_j$ with $z_{0j}$ for
$j=1,2$ and $y_1$ with $z_{10}$) and we can explicitly write down a global
section of $N_S$ in terms of an element of
$\Hom(\I_S/\I_S^2, \CO_S)$ as follows
\begin{equation}\label{D1ST}
\begin{split}
z_{i0} - y_1^i&\to 0\\
z_{i1} - x_1 y_1^i &\to -i h_2 y_1^{i-1}\\
z_{i2} - x_2 y_1^i &\to i h_1 y_1^{i-1}\\
f(x_1, x_2, y_1)&\to 0.
\end{split}
\end{equation}
The corresponding first order (embedded) deformation
$W_2\subset\P^{3k+2}\times \Spec\BC[t]/(t^2)$ of $S\subset \P^{3k+2}$ is 
defined locally on $Z_{00}\ne 0$ by
\begin{equation}\label{D1ST2}
\begin{split}
& z_{i0} = y_1^i\\
& z_{i1} = x_1 y_1^i + i h_2 y_1^{i-1} t\\
& z_{i2} = x_2 y_1^i - i h_1 y_1^{i-1} t\\
& f(x_1, x_2, y_1) = 0.
\end{split}
\end{equation}
To check \eqref{D1ST} defines a global morphism $\I_S/\I_S^2\to \CO_S$,
it is equivalent to check
\begin{quote}
($*$) the closure of the scheme defined by
\eqref{D1ST2} in $\P^{3k+2}\times \Spec\BC[t]/(t^2)$ has no component
other than $S$ as its central fiber over $\Spec\BC[t]/(t^2)$.
\end{quote}

First of all, $\P^2\times\P^1\subset\P^{3k+2}$ are covered by
affine open sets $Z_{01}\ne 0$, $Z_{02}\ne 0$, $Z_{k0}\ne 0$,
$Z_{k1}\ne 0$ and $Z_{k2}\ne 0$ besides $Z_{00}\ne 0$. And its local
defining functions over these affine sets are
\begin{gather*}
\frac{z_{ij} z_{01} - z_{i-1, j}z_{11}}{z_{01}^2},
\ \mathrm{for}\ i\ne 0\ \mathrm{and}\ (i,j)\ne (1,1)\ 
\mathrm{on}\ Z_{01}\ne 0\\
\frac{z_{ij} z_{02} - z_{i-1, j}z_{12}}{z_{02}^2}, \ \mathrm{for}\ i\ne 0\
\mathrm{and}\ (i,j)\ne (1,2)\
\mathrm{on}\ Z_{02}\ne 0\\
\frac{z_{ij} z_{k0} - z_{i+1, j}z_{k-1, 0}}{z_{k0}^2},\ \mathrm{for}\
i\ne k\ \mathrm{and}\ (i,j)\ne (k-1,0)\
\mathrm{on}\ Z_{k0}\ne 0\\
\frac{z_{ij} z_{k1} - z_{i+1, j}z_{k-1, 1}}{z_{k1}^2}, \text{ for }
i\ne k \text{ and } (i,j)\ne (k-1, 1)
\ \mathrm{on}\ Z_{k1}\ne 0\\
\frac{z_{ij} z_{k2} - z_{i+1, j}z_{k-1, 2}}{z_{k2}^2}, \text{ for }
i\ne k \text{ and } (i, j)\ne (k-1,2)
\ \mathrm{on}\ Z_{k2}\ne 0,
\end{gather*}
respectively. To verify ($*$), it suffices to take the closure of
\eqref{D1ST2} in these affine open sets one by one and see if the central
fiber consists of any component other than $S$. For example, on $Z_{01}\ne
0$, the closure of \eqref{D1ST2} is given by (plugging \eqref{D1ST2} into
$(z_{ij} z_{01} - z_{i-1, j}z_{11})/z_{01}^2$)
\begin{equation}\label{D1ST3}
\begin{split}
\frac{z_{i0} z_{01} - z_{i-1, 0}z_{11}}{z_{01}^2} &= -\frac{h_2
y_1^{i-1}}{x_1^2} t\\
\frac{z_{i1} z_{01} - z_{i-1, 1}z_{11}}{z_{01}^2} &= 0\\
\frac{z_{i2} z_{01} - z_{i-1, 1}z_{11}}{z_{01}^2} &= \frac{q
y_1^{i-1}}{x_1^2} t\\
\frac{f(x_1, x_2, y_1)}{x_1^3} &= 0
\end{split}
\end{equation}
where we identify $x_1, x_2, y_1$ with $z_{01}, z_{02}, z_{10}$ as before.
Since $h_2 y_1^{i-1}/x_1^2$ and $q y_1^{i-1}/x_1^2$ are regular over
$Z_{01}\ne 0$, the central fiber of \eqref{D1ST3} consists only of
$S$. The same analysis should be carried out on the other four open affine
sets for a complete verification of ($*$). But we will leave the details to
the readers.

The first order deformation $W_2$ of $S\subset\P^{3k+2}$ given in \eqref{D1ST2}
can be lifted to a deformation over disk $T$ since $H^1(N_S) = 0$
by \propref{P1}. Specifically, we can find $W\subset\P^{3k+2}\times
T$ with central fiber $S$ and locally defined by
\begin{align*}
z_{i0} &= y_1^i + O(t^2)\\
z_{i1} &= x_1 y_1^i + ih_2 y^{i-1} t + O(t^2)\\
z_{i2} &= x_2 y_1^i - ih_1 y^{i-1} t + O(t^2)
\end{align*}
and
$$f(x_1, x_2, y_1)= O(t^2).$$
Of course, there are infinitely many ways to lift a given first order
deformation. We only need to find one which serves our purpose.

\begin{claim}
We can inductively find $\alpha_{ij}(x_1, x_2, y_1)$,
$\beta_{ij}(x_1, x_2, y_1)$ and
$\gamma_{ij}(x_1, x_2, y_1)\in \BC[x_1, x_2, y_1]/(f(x_1, x_2, y_1))$
($1\le j\le i\le k$) with $\alpha_{i1} = 0$, $\beta_{i1} = i h_2$
and $\gamma_{i1} = -i h_1$
such that there exists $W\subset\P^{3k+2}\times T$ with central fiber $S$
and locally given by
\begin{equation}\label{DK3}
\begin{split}
z_{i0} &= y_1^i + \sum_{j=1}^i \alpha_{ij} y_1^{i-j+1} t^j + O(t^{i+1})\\
z_{i1} &= x_1 y_1^i + \sum_{j=1}^i \beta_{ij} y_1^{i-j}t^j + O(t^{i+1})\\
z_{i2} &= x_2 y_1^i + \sum_{j=1}^i \gamma_{ij} y_1^{i-j}t^j + O(t^{i+1})
\end{split}
\end{equation}
and
$$f(x_1, x_2, y_1) = 0.$$
\end{claim}

The procedure to find these polynomials can be described as follows.

Suppose we have found $\alpha_{ij}$, $\beta_{ij}$ and $\gamma_{ij}$ for
$j<l$. Namely, there exists $W_l\subset \P^{3k+2}\times \Spec\BC[t]/(t^l)$
locally defined by \eqref{DK3} up to order $t^{l-1}$.
We can lift $W_l$ to $W_{l+1}\subset \P^{3k+2}\times \Spec\BC[t]/(t^{l+1})$
locally defined by
\begin{equation}\label{DWL}
\begin{split}
z_{i0} &= y_1^i + \sum_{j=1}^{\min(i,l-1)} \alpha_{ij} y_1^{i-j+1} t^j +
\psi_i t^l + O(t^{i+2})\\
z_{i1} &= x_1 y_1^i + \sum_{j=1}^{\min(i,l-1)} \beta_{ij} y_1^{i-j}t^j 
+ \phi_i t^l + O(t^{i+2})\\
z_{i2} &= x_2 y_1^i + \sum_{j=1}^{\min(i,l-1)} \gamma_{ij} y_1^{i-j}t^j
+ \varphi_i t^l + O(t^{i+2})
\end{split}
\end{equation}
and
$$f(x_1, x_2, y_1) = f_l t^l$$
where $\psi_i, \phi_i, \varphi_i, f_l\in \BC[x_1, x_2, y_1]/(f(x_1, x_2, y_1))$
and we set $\psi_i = \phi_i = \varphi_i = 0$ for $i < l$. 
Since $W_{l+1}$ is an abitrary lift of $W_l$, we do not necessarily have
\begin{equation}\label{e:psi}
y_1^{i-l+1} | \psi_i, y_1^{i-l} | \phi_i, y_1^{i-l} | \varphi_i \
\text{and}\ f_l = 0.
\end{equation}
The idea here is to modify $\psi_i, \phi_i, \varphi_i, f_l$ one by one such
that \eqref{e:psi} holds and $W_{l+1}$ given locally by \eqref{DWL}
remains as a deformation of $S$ over $\Spec\BC[t]/(t^{l+1})$.

Again, to check that $W_{l+1}$ given locally by \eqref{DWL} is a flat
family of surfaces over $\Spec\BC[t]/(t^{l+1})$ with central fiber $S$,
it is equivalent to check that
the closure of the scheme defined by
\eqref{DWL} in $\P^{3k+2}\times \Spec\BC[t]/(t^{l+1})$ has no component
other than $S$ as its central fiber.

First, we can obviously set $f_l = 0$ and
hence inductively we can set $f(x_1, x_2,
y_1) = 0$ in the lift of $W_2$ of any order (this also follows from
$H^1(N_{\P^2\times \P^1}|_S) = 0$ as proved in \propref{P1}).

Take any polynomial $g$ lying the $\BC$-linear space
spanned by $y_1^i$, $x_1 y_1^i$ and $x_2 y_1^i$ ($0\le i\le k$). Our
first observation is that we may replace any $\psi_i$ ($\phi_i$ or
$\varphi_i$) for $i\ge l$ by $\psi_i+g$ ($\phi_i+g$ or
$\varphi_i+g$) and the corresponding
$W_{l+1}$ locally defined by \eqref{DWL} is still a lift of
$W_l$. Therefore, it suffices to show that for each $\psi_i$ ($\phi_i$ or
$\varphi_i$) 
there exists $g\in\oplus_{j=0}^k
(\BC y_1^j \oplus \BC x_1 y_1^j\oplus \BC x_2 y_1^j)$ such that
$y_1^{i-l+1} | (\psi_i + g)$ ($y_1^{i-l} | (\phi_i + g)$ or $y_1^{i-l} |
(\varphi_i + g)$).

The defining equations of $W_{l+1}$ over $Z_{01}\ne 0$ can be obtained by
taking the closure of \eqref{DWL} over $Z_{01}\ne 0$. As it is illustrated
in the case $l=1$, this is done by simply plugging \eqref{DWL} into $(z_{ij}
z_{01} - z_{i-1, j}z_{11})/z_{01}^2$. We will get a set of equations in the
form
\begin{equation*}
\begin{split}
\frac{z_{i0} z_{01} - z_{i-1, 0}z_{11}}{z_{01}^2} &= \sum_{j=1}^l a_{ij} t^j\\
\frac{z_{i1} z_{01} - z_{i-1, 1}z_{11}}{z_{01}^2} &= \sum_{j=1}^l b_{ij} t^j\\
\frac{z_{i2} z_{01} - z_{i-1, 1}z_{11}}{z_{01}^2} &= \sum_{j=1}^l c_{ij} t^j.
\end{split}
\end{equation*}
Our calculation shows that
\begin{align*}
a_{il} &= \frac{x_1 \psi_i - x_1 \psi_{i-1} y_1 - h_2 \alpha_{i,l-1}
y_1^{i-l+1} + O(y_1^{i-l+1})}{x_1^2}\\
b_{il} &= \frac{x_1 \phi_i - x_1 \phi_{i-1} y_1 - h_2 \beta_{i,l-1}
y_1^{i-l} + O(y_1^{i-l+1})}{x_1^2}\\
c_{il} &= \frac{x_1 \varphi_i - x_1 \varphi_{i-1} y_1 - h_2 \gamma_{i,l-1}
y_1^{i-l} + O(y_1^{i-l+1})}{x_1^2}.
\end{align*}
Since the central fiber of $W_{l+1}$ consists only of $S$,
$a_{il}$ must be a regular function over $Z_{01}\ne 0$. And since
$\BC[x_1, x_2, y_1]/(f(x_1, x_2, y_1))$ is a UFD, we must have
$$x_1 \psi_i - x_1 \psi_{i-1} y_1 - h_2 \alpha_{i,l-1}
y_1^{i-l+1} + O(y_1^{i-l+2}) = p(x_1, x_2, y_1)$$
in $\BC[x_1, x_2, y_1]/(f(x_1, x_2, y_1))$
for some polynomial $p(x_1, x_2, y_1)$ which is quadratic in $x_1$ and
$x_2$. By induction on $i$, we may assume $y_1^{i-l} | \psi_{i-1}$. Let
$m$ be the number such that $y_1^m | \psi_i$ and $y_1^{m+1} \nmid
\psi_i$.
If $m \ge i-l+1$, we are done. If not, obviously we have $y_1^m |
p(x_1, x_2, y_1)$.
Let $\psi_i(x_1, x_2, y_1)
= y_1^m \delta(x_1, x_2, y_1)$ and $p(x_1, x_2, y_1) = y_1^m r(x_1, x_2,
y_1)$. Then
$$x_1\delta(x_1, x_2, 0) = r(x_1, x_2, 0)$$
in the ring $\BC[x_1, x_2]/(f(x_1, x_2, 0))$. Since we assume the point
$(X_0 = X_1 = Y_1 = 0)\not\in S$ and $r(x_1, x_2, 0)$ is a quadratic
polynomial in $x_1$ and $x_2$, we must have
$$r(x_1, x_2, 0) = x_1 \lambda(x_1, x_2)$$
for some $\lambda(x_1, x_2)\in\BC\oplus\BC x_1\oplus \BC x_2$. So we may
replace $\psi_i$ by $\psi_i - \lambda(x_1, x_2)y_1^m$ which is easy to see
divisible by $y^{m+1}$.
We can repeat this procedure until $y_1^{i-l+1} | \psi_i$. And
we can do the same to $\phi_i$ and $\varphi_i$. Hence we eventually arrive
at $W\subset \P^{3k+2}\times T$ which is locally given by \eqref{DK3}.

Furthermore, if we let $\{\lambda_i\}$ be a sequence of elements in the ring
$\BC[x_1, x_2]/(f(x_1, x_2, 0))$ satisfying the recursive condition
\begin{subequations}
\begin{equation}\label{E:LAMBDA}
\lambda_1 = h_2(x_1, x_2, 0) \ \text{and}\ x_1\lambda_{i+1} = h_2(x_1, x_2,
0)\lambda_i + r_i(x_1, x_2)
\end{equation}
where $r_i(x_1, x_2)$ is some quadratic polynomial in $x_1$ and $x_2$,
we may choose $\beta_{ij}(x_1, x_2, y_1)$ such that $\beta_{ij}(x_1, x_2,
0) = \binom{i}{j} \lambda_j$ by the same argument as above. Similarly, 
if we let $\{\mu_i\}$ be a sequence of elements in the ring
$\BC[x_1, x_2]/(f(x_1, x_2, 0))$ satisfying the recursive condition
\begin{equation}\label{E:MU}
\mu_1 = -h_1(x_1, x_2, 0) \ \text{and}\ x_1\mu_{i+1} = h_2(x_1, x_2,
0)\mu_i + s_i(x_1, x_2)
\end{equation}
where $s_i(x_1, x_2)$ is some quadratic polynomial in $x_1$ and $x_2$,
we may choose $\gamma_{ij}(x_1, x_2, y_1)$ such that $\gamma_{ij}(x_1, x_2,
0) = \binom{i}{j} \mu_j$. Hence there exists a flat family $W\subset
\P^{3k+2}\times T$ over $T$ which is locally given by
\begin{equation}\label{E:DK3}
\begin{split}
& z_{i0} = y_1^i + \sum_{j=1}^i O(y_1^{i-j+1}) t^j + O(t^{i+1})\\
& z_{i1} = x_1 y_1^i + \sum_{j=1}^i \left(\binom{i}{j}\lambda_j +
O(y_1)\right) y_1^{i-j}t^j + O(t^{i+1})\\
& z_{i2} = x_2 y_1^i + \sum_{j=1}^i \left(\binom{i}{j}\mu_j + O(y_1)\right)
y_1^{i-j}t^j + O(t^{i+1})\\
& f(x_1, x_2, y_1) = 0.
\end{split}
\end{equation}
\end{subequations}
And it is not hard to see that the general fibers of $W$ are primitive K3
surfaces in $\P^{3k+2}$.

Now we are ready to prove \thmref{T:NODB}.

Let $\Upsilon_t$ be a general fiber of $\Upsilon \to T$.
Let $\delta(\Upsilon_t, Z)$ be the total $\delta$-invariant of $\Upsilon_t$
in the neighborhood of $Z\subset\Upsilon_0$.
We will show that 
\begin{claim}\label{clm:2}
$\delta(\Upsilon_t, \Gamma_i) > 3m_i$ ($i>0$) for $W$ general.
\end{claim}
Notice that if this is true for some $W$ whose general fibers are
primitive K3 surfaces, it should be true for $W$ general.

It is obvious that $\delta(\Upsilon_t, \Gamma_i) = 3$ if $m_i = 1$.
And since $\delta(\Upsilon_t, (\Gamma_0)_{sing}) = 3l+2$ where
$(\Gamma_0)_{sing}$ is the singular locus of $\Gamma_0$, the total
$\delta$-invariant of $\Upsilon_t$ will exceed $3k+2$ if $m_i > 1$ for some
$i$, namely, $\Gamma$ is nonreduced. This is impossible since we know
$\Upsilon_t$ is an irreducible rational curve with arithemtic genus $3k+2$.
Hence \claimref{clm:2} directly implies that $\Gamma$ is reduced. 

So let us assume $m = m_1 > 1$.
Let $\Gamma_0\cap\Gamma_1 = \{P_1\} \cup \{P_2\} \cup \{P_3\}$
and $Q$ be the node of $\Gamma_1$.

After an action of $\P GL(3)\times\P GL(2)$ on $S\subset\P^2\times\P^1$
which is induced by $\P GL(3k+3, 3k+3)$ on $S\subset\P^{3k+3}$, we may assume that $\Gamma_0$ is
given by $x_1 + O(y_1) = 0$
and $\Gamma_1$ is given by $y_1 = 0$ on $S$. And by an action of
$\P GL(3)\times\P GL(2)$ fixing $x_1 = 0$ and $y_1=0$, we can make $S$ miss
points $(X_0=X_1=Y_1=0)$ and $(X_0=X_2=Y_1=0)$ as required
before.

Let $W$ be the family of K3 surfaces locally given by \eqref{E:DK3} up to base changes.
We may make the following assumptions (the reason we do
so will be clear in a moment)
\begin{enumerate}
\item The line $X_0 = 0$ on the plane $Y_1 = 0$
meets $\Gamma_1$ at three distinct points
$R_1$, $R_2$ and $R_3$.
\item Let
$$\kappa = \frac{h_2(x_1, x_2, 0)}{x_1^2} \ \text{and}\ \rho =
\frac{r_1(x_1, x_2)}{x_1^2}$$
for $r_1(x_1, x_2)$ in \eqref{E:LAMBDA}.
Notice that we have some freedom to choose $h_1$ and $h_2$ in
\eqref{E32}. We may replace $h_1$ and $h_2$ by $h_1+x_2l(x_1, x_2)$ and
$h_2 - x_1 l(x_1, x_2)$ for any $l(x_1, x_2)\in \BC\oplus \BC x_1\oplus \BC
x_2$. Similarly, we have some freedom to choose $\{\lambda_i\}$ and
$\{r_i(x_1,x_2)\}$ in \eqref{E:LAMBDA}. We may replace $r_1(x_1, x_2)$ by
$r_1(x_1, x_2) + x_1 l(x_1, x_2)$ for any $l(x_1, x_2)\in \BC\oplus \BC
x_1\oplus \BC x_2$ (and change $\lambda_i$ for $i>1$ accordingly).
Hence we can make $\kappa$ not vanish at $R_1$, $R_2$ and $R_3$ and
the values of $\rho \kappa^{-1}$ at $R_1$, $R_2$ and $R_3$
different from each other.
\end{enumerate}
Let $\Upsilon$ be cut out by a family of
hyperplanes $H_t$ in $\P^{3k+2}$,
which is given by (after a proper base change)
\begin{equation}\label{E:HT}
z_{m1} + \sum_{i,j} a_{ij}(t) z_{ij} = 0
\end{equation}
where $a_{ij}(t)\in \BC[[t]]$ and $a_{ij}(0) = 0$ for $i \le m$.
Combining
\eqref{E:DK3} and \eqref{E:HT} (notice we have made a base change so
$t$ in \eqref{DK3} should be replaced by $t^l$ for some $l$), we have the
defining equation of $\Upsilon_t$ on $Z_{00}\ne 0$
\begin{align}\label{E:UPS}
&\quad\left(x_1y_1^m+\sum_{j=1}^m
\beta_{mj}y_1^{m-j}t^{lj}+O(t^{l(m+1)})\right)\\
&+\sum a_{i0}(t)
\left (y_1^i + \sum_{j=1}^i O(y_1^{i-j+1}) t^{lj} 
+ O(t^{l(i+1)})\right)\notag\\
&+\sum a_{i1}(t)
\left (x_1 y_1^i + \sum_{j=1}^i \beta_{ij}
y_1^{i-j} t^{lj} 
+ O(t^{l(i+1)})\right)\notag\\
&+\sum a_{i2}(t)
\left (x_2 y_1^i + \sum_{j=1}^i \gamma_{ij}
y_1^{i-j} t^{lj} 
+ O(t^{l(i+1)})\right)=0,\notag
\end{align}
where
$$\beta_{ij} = \binom{i}{j}\lambda_j + O(y_1) \ \text{and}\
\gamma_{ij} = \binom{i}{j}\mu_j + O(y_1).$$
Let
$$\delta = \min\left(l,
\min_{i<m}\left(\frac{\nu(a_{ij}(t))}{m-i}\right)\right)\ \text{and}\ c_{ij}
= \left.\frac{a_{ij}(t)}{t^{(m-i)\delta}}\right|_{t=0}$$
where $\nu(a_{ij}(t))$ is the valuation of $a_{ij}(t)\in \BC[[t]]$.
We may make $\delta$ an integer by a proper base change. Let $\pi:
\Upsilon'\to\Upsilon$
be the blow-up of $\Upsilon$ along the subscheme $y_1 = t^\delta = 0$
and let $y = y_1 / t^\delta$, 
$\Gamma_0' = \overline{\pi^{-1}(\Gamma_0-\{P_1, P_2, P_3\})}$ and
$\Gamma_1' = \pi^{-1}(\Gamma_1) \subset \Upsilon_0'$.
And let $\Upsilon^v$ be the nodal reduction of $\Upsilon'$.

The curve $\Gamma_0'$ can be described as a curve in
$\P^2\times\P^1$ with affine coordinates $(x_1, x_2)\times (y)$.
There are two cases to consider.

\subsection{The case $\delta < l$} 
By \eqref{E:UPS},
$\Gamma_1'$ is given by
$$
x_1 y^m + l_1(x_1, x_2) y^{m-1} + l_2(x_1, x_2) y^{m-2} + ... + l_r(x_1,
x_2) y^{m-r} = 0
$$
and
$$f(x_1, x_2, 0) = 0$$
where $l_i(x_1, x_2)= c_{i0} + c_{i1}x_1 + c_{i2}x_2$
for $i=1,
2,...,r$. Obviously, $\Gamma_1' = C_1\cup (m-r)C_2$ where $C_1$ is a
curve given by
\begin{equation}\label{E:C1}
x_1y^r + \sum_{i=1}^r l_i(x_1, x_2) y^{r-i} = 0
\end{equation}
which maps to
$\Gamma_1$ with degree $r$ and $C_2$ is given by $y=0$.

\subsubsection{}
If $l_r(x_1, x_2) \not\in \BC x_1$, then the line $l_r(x_1, x_2) = 0$
in the plane $y_1 = 0$ passes through
at most one of the three points $P_1$, $P_2$
and $P_3$, say it misses $P_1$ and $P_2$. Factoring the LHS of
\eqref{E:C1} as a polynomial in $y$
over the ring $\BC[x_1,x_2]/(f(x_1,x_2,0))$, we
get an irreducible component $C_1'$ of $C_1$ given by
$$x_1y^{s} + \sum_{i=1}^s l_i'(x_1, x_2) y^{s-i} = 0.$$
Since $l_r(x_1, x_2) \ne 0$ at $P_1$ and $P_2$, $l_s'(x_1, x_2) \ne 0$ at
$P_1$ and $P_2$, either. Hence $C_1'$ meets $\Gamma_0'$ at two
points $P_1'\in \pi^{-1}(P_1)$ and $P_2'\in \pi^{-1}(P_2)$
where $u = 1/y = t^{\delta}/y_1 = 0$ and $\Upsilon$ is locally given by
\begin{gather*}
x_1 + \sum_{i=1}^s l_i'(x_1, x_2) u^i = O(t)\\
u y_1 = t^\delta
\end{gather*}
and
$$f(x_1, x_2, y_1) = 0.$$
Obviously, the general fiber $\Upsilon_t'$ is smooth in the neighborhoods of
$P_1'$ and $P_2'$.
This implies that $\wt{C_1'}$ will have at least two
intersections with $\wt{\Gamma_0'}$, where $\wt{C_1'}$ and $\wt{\Gamma_0'}$
are the irreducible components of the central fiber $\Upsilon_0^v$ of
$\Upsilon^v$ which dominates $C_1'$ and
$\Gamma_1'$, respectively. This is impossible since the general fiber
$\Upsilon^v$ is rational.

\subsubsection{}
If $l_r(x_1, x_2)\in \BC x_1$, $l_r(x_1, x_2)$ does not vanish at the
node $Q$ of $\Gamma_1$. Choose point $Q'\in \pi^{-1}(Q)$ and let
$C_1'$ be the union of components of $C_1$ passing through
$Q'$. Let $u$ and $v$ be the
local coordinates of $\Gamma_1$ at $Q$ such that $u=v=0$ at $Q$ and
$f(x_1, x_2, 0) = uv$. Let
$a$ be the $y$-coordinate of point $Q'$. Since $l_r(x_1, x_2)\ne 0$ at $Q$,
$a\ne 0$. Hence in a neighborhood of $Q'$, $\Upsilon$ is locally given by
\begin{equation}\label{E:QP}
\begin{split}
(y-a)^s &= O(u, v, t)\\
uv &= t^\delta(1+O(u,v,y-a)) + O(t^{2\delta}),
\end{split}
\end{equation}
where $s$ is ramification index of $Q'$ under $\pi$.
Let $\Sigma$
be the union of irreducible components of $\Upsilon_0^v$
which map nonconstantly to $C_1'$ and hence dominate $\Gamma_1$. The
morphism $\Sigma\to \Gamma_1$ must factor through $\wt{\Gamma_1}$, where
$\wt{\Gamma_1}$ be the normalization of $\Gamma_1$. Let $\phi$ be the
morphism $\Sigma\to\wt{\Gamma_1}$ and $\varphi$ be the morphism $\Sigma\to
C_1'$. And let
$Q_1$ and $Q_2$ be two points on $\wt{\Gamma_1}$ mapping to $Q$
corresponding to the branches $u = 0$ and $v = 0$, respectively, and 
$$S_1 = \phi^{-1}(Q_1)\cap \varphi^{-1}(Q')\ \text{and}\ S_2 =
\phi^{-1}(Q_2)\cap \varphi^{-1}(Q').$$
By \eqref{E:QP} and \lemref{L:XY=T}, each point $p\in S_1$ ($S_2$)
is joined by a chain of curves not in $\Sigma$
to some point $q\in S_2$ ($S_1$).
Since each irreducible component of $\Sigma$ has points
in both $S_1$ and $S_2$, the components of $\Sigma$
forms a graph in which each vertex (representing a component)
has degree at least two and hence the graph must contain a
cycle. This contradicts the fact that the components of $\Upsilon_0^v$
forms a tree.

\subsection{The case $\delta\ge l$} By \eqref{E:UPS}
$\Gamma_1'$ is given by
\begin{align}\label{E:GAMMA1P}
&\quad\left(x_1 y^m + \sum_{j=1}^m \binom{m}{j} \lambda_j y^{m-j}\right) +
\sum_{i>0} c_{i0} y^{m-i}\\
& +\sum_{i>0} c_{i1} \left (x_1 y^{m-i} 
+ \sum_{j=1}^{m-i} \binom{m-i}{j} \lambda_j y^{m-i-j}\right)\notag\\
& +\sum_{i>0} c_{i2} \left (x_2 y^{m-i} 
+ \sum_{j=1}^{m-i} \binom{m-i}{j} \mu_j y^{m-i-j}\right) = 0\notag
\end{align}
and
\begin{equation}\label{E:GAMMA1P2}
f(x_1, x_2, 0) = 0.
\end{equation}
Let $\Gamma_1' = Z\cup \Gamma_1''$ where $Z$ is the union of components of
$\Gamma_1'$ which map constantly to the points $P_i$ for $i=1,2,3$. Obviously,
$Z$ is reduced and consists of at most three lines corresponding to $P_i$.
We claim that
\begin{claim}\label{clm:gamma1p} $\Gamma_1''$ is reduced for any $m\ge 2$.
And $\Gamma_1''$ is either irreducible or
consisting of two irreducible
components, each mapping to $\Gamma_1$ with degree $m/2$ ($m$ must be even
in this case). 
\end{claim}
This is done by showing that \eqref{E:GAMMA1P} as a polynomial in $y$ over the
function field $K(\Gamma_1)$ of $\Gamma_1$ is either irreducible
or factored into two distinct polynomials, each with degree $m/2$. While
this is in turn proved by
localizing $\Gamma_1$ at one of the points $R_1$, $R_2$ and
$R_3$. Since $\rho\kappa^{-1}$ has different values at $R_1$, $R_2$ and
$R_3$, we may assume that
\begin{subequations}
\begin{gather}\label{E:A1}
c_{10} - (m-1)\rho\kappa^{-1}\ne 0\\
\intertext{and}
\label{E:A2}
(c_{11} + c_{12}\mu)^2 - 4(c_{11} + c_{12}\mu + \rho -
c_{10}\kappa) \ne 0
\end{gather}
\end{subequations}
hold simultaneously at one of $R_i$, say $R_1$,
where $\kappa = h_2(x_1, x_2, 0) / x_1^2$ and $\rho = r_1(x_1, x_2) /
x_1^2$ as defined before.
Localize $\Gamma_1''$ at $R_1$ and let $\CO_{\Gamma_1, R_1}$ be the
local ring of $\Gamma_1$ at $R_1$ with uniformizer
$\varepsilon = 1/x_1$ and let $\eta = x_2/x_1$.
According to \eqref{E:LAMBDA} and \eqref{E:MU}, we have
\begin{equation*}
\begin{split}
& \frac{\lambda_1}{x_1^2} = \kappa,\ \text{and}\
\frac{\lambda_{i+1}}{x_1^{i+2}} = (\kappa^2 +
\rho \varepsilon^2)\kappa^{i-1} + O(\varepsilon^3), \
\text{for}\ i\ge 1\\
& \frac{\mu_i}{x_1^{i+1}} = \eta\kappa^i + O(\varepsilon^2).
\end{split}
\end{equation*}
Hence
\begin{align*}
& \frac{1}{x_1^{m+1}} \left (x_1 y^m  + \sum_{j=1}^m \binom{m}{j} \lambda_j
y^{m-j}\right)\\
& \quad\quad = (1+\rho\kappa^{-2}\varepsilon^2) (\varepsilon y + \kappa)^m
- \rho\kappa^{-2}\varepsilon^2\left((\varepsilon y)^m + m \kappa (\varepsilon
y)^{m-1}\right)\\
& \quad\quad\quad+ \sum_{i>1} O(\varepsilon^3) (\varepsilon y)^{m-i},\\
& \frac{1}{x_1^{m+1}} \sum_{i>0} c_{i0} y^{m-i} = c_{10} \varepsilon^2
(\varepsilon y)^{m-1} + \sum_{i>1} O(\varepsilon^3) (\varepsilon
y)^{m-i},\\
& \frac{1}{x_1^{m+1}} \sum_{i>0} c_{i1} \left (x_1 y^{m-i} 
+ \sum_{j=1}^{m-i} \binom{m-i}{j} \lambda_j y^{m-i-j}\right)\\
& \quad\quad = c_{11}\varepsilon (\varepsilon y + \kappa)^{m-1}
+ c_{21} \varepsilon^2 (\varepsilon y + \kappa)^{m-2}
+\sum_{i>1} O(\varepsilon^3) (\varepsilon y)^{m-i},\\
& \frac{1}{x_1^{m+1}} \sum_{i>0} c_{i2} \left (x_2 y^{m-i} 
+ \sum_{j=1}^{m-i} \binom{m-i}{j} \mu_j y^{m-i-j}\right)\\
&\quad\quad = c_{12}\eta\varepsilon (\varepsilon y + \kappa)^{m-1}
+ c_{22} \eta\varepsilon^2 (\varepsilon y + \kappa)^{m-2}
+\sum_{i>1} O(\varepsilon^3) (\varepsilon y)^{m-i}.
\end{align*}
Adding these equations up, we obtain the defining equation of $\Gamma_1''$
as a scheme over the local ring $\CO_{\Gamma_1, R_1}$ with uniformizer
$\varepsilon$
\begin{equation}\label{E:GAMMA1PL}
\begin{split}
w^m & + (c_{11}+c_{12}\eta)\varepsilon
w^{m-1} + (c_{21}+c_{22}\eta)\varepsilon^2 w^{m-2}\\
& + \varepsilon^2 \left((c_{10} - m\rho\kappa^{-1}) 
(w-\kappa)^{m-1} - \rho\kappa^{-2}(w-\kappa)^m\right) \\
& \quad + \sum_{i=0}^m O(\varepsilon^3) w^i = 0
\end{split} 
\end{equation}
where $w = \varepsilon y + \kappa$. Since $(c_{10} -
m\rho\kappa^{-1})(-\kappa)^{m-1} - \rho\kappa^{-2}(-\kappa)^m\ne 0$
valuated at $R_1$ by \eqref{E:A1},
the LHS of \eqref{E:GAMMA1PL}, as a polynomial in $w$,
has $m$ different roots in $\overline{K(\Gamma_1)}$, each at the
order of $\varepsilon^{2/m}$ for $m>2$. So $\Gamma_1''$ is reduced for
$m>2$. Furthermore,
it is not hard to see that the LHS of \eqref{E:GAMMA1PL} is an
irreducible polynomial in $w$ for $m$ odd and it is either irreducible or
factored into two irreducible
polynomials of degree $m/2$ each for $m$ even and $m>2$. It remains to show
that $\Gamma_1''$ is reduced for $m=2$. Taking the discriminant
of the LHS of \eqref{E:GAMMA1PL} for $m=2$, we have
$$\left((c_{11} + c_{12}\mu)^2 - 4(c_{11} + c_{12}\mu + \rho -
c_{10}\kappa)\right)\varepsilon^2 + O(\varepsilon^3)$$
which is nonzero by \eqref{E:A2}. Hence
$\Gamma_1''$ is reduced for $m=2$. We have established \claimref{clm:gamma1p}.

Therefore, $\Gamma_1'$ is a reduced curve in $\P^2\times\P^1$ which is the
complete intersection of two surfaces \eqref{E:GAMMA1P} and
\eqref{E:GAMMA1P2} of type $(m+1, m)$ and $(3,0)$, respectively. Hence
$p_a(\Gamma_1') = 3m^2 - 2$. And since $\Gamma_0'$ meets $\Gamma_1'$ at
three points, $\delta(\Upsilon_t', \Gamma_1') = 3m^2 > 3m$ when $m > 1$.
This finishes the proof of \claimref{clm:2} and hence \thmref{T:NODB}.

\begin{cor}\label{COR:NODB}
\conjref{T:NODA} holds true if $l\le 2$ for {\bf{TK1}}, $l\le 1$ for
{\bf{TK2}}, or $l = 0$ for {\bf{TK3}} and \conjref{T:NODC} holds true if
$l\le 3$ for {\bf{TK1}}, $l\le 2$ for {\bf{TK2}} or $l \le 1$ for {\bf{TK3}}.
Hence by \thmref{T:NODB}, \conjref{T:NOD} holds true
for $n\le 9$ and $n=11$.
\end{cor}

\begin{proof}
The argument is similar to that used in the proof of \conjref{T:NOD} for the
cases that $n=3, 4$.

Let $W\subset 
|C + lF| \times |-K_{\P E}|$ be the incidence
correspondence $(H, S)$ such that 
$H\cap S$ is an irreducible rational curve, where $K_{\P
E}$ is the canonical divisor of $\P E$. It is not
hard to see that $H^1(-K_{\P E} - C - lF) = 0$ if $l\le 3$ for {\bf{TK1}},
$l\le 2$ for {\bf{TK2}} or $l \le 1$ for {\bf{TK3}}. Hence $|-K_{\P E}|$
cuts out the complete linear series $|\CO_H(-K_{\P E})|$ on $H\in
|C+lF|$. Therefore, the fiber $W_H$ of $W$ over $H\in |C + lF|$ can be
identified with $V\times 
H^0(-K_{\P E} - C - lF)$, where $V$ is the variety parametrizing irreducible
ratonal curves in $|\CO_H(-K_{\P E})|$. By \thmref{T:HARRIS}, the general
member of $V$ is nodal. Hence \conjref{T:NODC} holds true if $l\le 3$ for
{\bf{TK1}}, $l\le 2$ for {\bf{TK2}} or $l \le 1$ for {\bf{TK3}}.

Besides, since $\dim V > 1$, a general member $C\in V$ meets a curve in
$|\CO_H(F)|$ transversely at three points. Let us fix $C\in V$ and $D\in
|F|$ such that $C$ meets $D\isom \P^2$ at three distinct points
$p_1$, $p_2$ and
$p_3$. Considering $(H, S)\in \{C\}\times H^0(-K_{\P E} - C - lF)\subset
W_H$, we see that such $S$ cuts out a linear series 
$\sigma\subset |\CO_D(-K_{\P
E})| = |\CO_{\P^2}(3)|$ with base points at $p_1$, $p_2$ and $p_3$ if
$l\le 2$ for {\bf{TK1}}, $l\le 1$ for {\bf{TK2}}, or $l = 0$ for
{\bf{TK3}}. Obviously, a general irreducible rational curve in $\sigma$ is also
a general irreducible rational curve in $|\CO_{\P^2}(3)|$. Hence
\conjref{T:NODA} holds if $l\le 2$ for {\bf{TK1}}, $l\le 1$ for {\bf{TK2}},
or $l = 0$ for {\bf{TK3}}.
\end{proof}

\end{document}